\documentclass[a4paper,10pt]{article}
\usepackage{a4, amsxtra,amsmath,amsfonts,amscd, amssymb, mathrsfs}
\usepackage[all]{xy}
\newtheorem{thm}{Theorem}[section]
\newtheorem{prop}[thm]{Proposition}
\newtheorem{cor}[thm]{Corollary}
\newtheorem{lem}[thm]{Lemma}
\newtheorem{ex}[thm]{Example}
\newtheorem{Def}[thm]{Definition}
\newtheorem{rem}[thm]{Remark}
\newtheorem{art}[thm]{}

\newcommand{\Div}{{\rm div}}

\newcommand{\Pic}{{\rm Pic}}
\newcommand{\ord}{{\rm ord}}
\newcommand{\Spec}{{\rm Spec}}

\newcommand{\ve}{{\varepsilon}}

\newcommand{\Gal}{{\rm Gal}}
\newcommand{\supp}{{\rm supp}}

\newcommand{\Hcal}{{\mathscr H}}
\newcommand{\Jcal}{{\mathscr J}}
\newcommand{\Lcal}{{\mathscr L}}
\newcommand{\Mcal}{{\mathscr M}}
\newcommand{\Ncal}{{\mathscr N}}
\newcommand{\Ocal}{{\mathscr O}}

\newcommand{\Vcal}{{\mathscr V}}
\newcommand{\Xcal}{{\mathscr X}}
\newcommand{\Ycal}{{\mathscr Y}}
\newcommand{\cdop}{{\mathbb C}}

\newcommand{\qdop}{{\mathbb Q}}
\newcommand{\ndop}{{\mathbb N}}
\newcommand{\rdop}{{\mathbb R}}
\newcommand{\pdop}{{\mathbb P}}
\newcommand{\kdop}{{{\mathbb K}_v}}

\newcommand{\zdop}{{\mathbb Z}}

\newcommand{\metr}{{\|\hspace{1ex}\|}}
\newcommand{\canmetr}{{\|\hspace{1ex}\|_{\rm can}}}
\newcommand{\metrv}{{\|\hspace{1ex}\|_v}}
\newcommand{\absv}{{|\hspace{1ex}|_v}}

\newcommand{\proof}{\noindent {\bf Proof: \/}}
\newcommand{\qed}{{ \hfill $\square$}}

\newcommand{\val}{{\rm val}}

\newcommand{\xb}{{\mathbf x}}

\newcommand{\Tor}{{\mathbb G}_m^n}
\newcommand{\rtor}{{\rdop^n/\Lambda}}

\newcommand{\Deltabar}{{\overline{\Delta}}}

\newcommand{\valbar}{{\overline{\val}}}
\newcommand{\Xan}{{X_v^{\rm an}}}
\newcommand{\Zan}{{Z_v^{\rm an}}}

\newcommand{\kcirc}{{{\mathbb K}_v^\circ}}

\title{Equidistribution over function fields}
\author{Walter Gubler}
\date{\today}
\begin{document}

\maketitle

\begin{abstract}
We prove equidistribution of a generic net of small points in a projective variety $X$ over a function field $K$. For an algebraic dynamical system over $K$, we generalize this equidistribution theorem to a small generic net of subvarieties. For number fields, these results were proved by Yuan and we transfer here his methods to  function fields. If $X$ is a closed subvariety of an abelian variety, then we can describe the equidistribution measure explicitly in terms of convex geometry.
\end{abstract} 

\section{Introduction}

Equidistribution plays an important role in various branches of science. In diophantine geometry, equidistribution of small points is the key in Ullmo's and Zhang's proof of the Bogomolov conjecture (see \cite{Ullm}, \cite{Zh}). The goal of this article is to transfer Yuan's recent generalization of this equidistribution theorem to function fields. Since the Bogomolov conjecture over function fields is still open, this might have interesting applications.

In this paper, we consider the function field $K=k(B)$ of an integral projective variety $B$ over the field $k$ with $B$  regular in codimension $1$.  We fix an ample class $\mathbf c$ on $B$ to get an absolute height function $h_{(L,\metr)}$ on an irreducible $d$-dimensional projective variety $X$ over $K$ with respect to an admissibly metrized line bundle $(L,\metr)$ (see \S 3 for details). We fix also a place $v$ of $K$. The Berkovich analytic space $\Xan$ is defined over the smallest algebraically closed field $\kdop$ which is complete with respect to a valuation extending $v$. There is a regular Borel measure $c_1(L,\metr_v)^{\wedge d}$ on $\Xan$ which was introduced by Chambert-Loir and which is the non-archimedean analogue of the corresponding complex analytic wedge product of Chern forms (see \S 2).

We choose always an embedding $\overline{K} \hookrightarrow \kdop$ to identify $X(\overline{K})$ with a subset of $\Xan$. The Galois group $\Gal(\overline K/K)$ acts on $X$ and we denote the orbit of $P \in X(\overline K)$ by $O(P)$. Let $\delta_P$ be the Dirac measure on $\Xan$ in $P$. We consider a {\it generic net} $(P_m)_{m \in I}$ in $X(\overline K)$, i.e. $I$ is a directed set and for every proper closed subset $Y$ of $X$, there is $m_0 \in I$ with $P_m \not \in Y$ for all $m \geq m_0$. Our main result is the following {\it equidistribution theorem}:

\begin{thm} \label{main theorem}
Let $L$ be a big semiample line bundle on the irreducible $d$-dimen\-sional projective variety $X$ over the function field $K$. We endow $L$ with a semipositive admissible metric $\metr$. We assume that $(P_m)_{m \in I}$ is a generic net in $X(\overline K)$ with 
$$\lim_m h_{(L,\metr)}(P_m)=\frac{1}{(d+1)\deg_L(X)} h_{(L,\metr)}(X).$$
For a place $v$ of $K$, we have the following weak limit of regular probability measures on $\Xan$:
$$\frac{1}{|O(P_m)|} \sum_{P_m^\sigma \in O(P_m)} \delta_{P_m^\sigma} \stackrel{w}{\rightarrow} \frac{1}{\deg_L(X)} c_1(L,\metr_v)^{\wedge d}.$$
\end{thm}

For number fields, the history of equidistribution of small points starts with the theorem of Szpiro--Ullmo-Zhang \cite{SUZ} which handles the case of an archimedean place $v$ assuming positive curvature and $X$ smooth. Smoothness was removed by Ullmo \cite{Ullm} and Zhang \cite{Zh} to prove the Bogomolov conjecture. Chambert-Loir has proved a non-archimedean version of these results assuming that the metric at the finite place $v$ of the number field is induced by an ample model (see \cite{Ch}, Th\'eor\`eme  3.1). In case of abelian varieties and canoncial metrics, this handles just the case of good reduction at $v$ as the canonical metric is only semipositive at places of bad reduction. Moreover, Chambert-Loir  proved  equidistribution of small points for any admissible metric in the case of curves (\cite{Ch}, Th\'eor\`eme 4.2). Finally, Yuan introduced a variational principle  to deal also with semipositive metrics and he proved Theorem \ref{main theorem} for any place $v$ of the number field $K$ (see \cite{Yu}, \S 5).

For function fields, a tropical version of Theorem \ref{main theorem} was proved in case of a closed subvariety of an abelian variety which is totally degenerate at the given place $v$ (see \cite{Gu5}, Theorem 5.5). This was the key to prove Bogomolov's conjecture for totally degenerate abelian varieties over function fields. Independently of the present paper, Faber (\cite{Fa}, Theorem 1.1) proved equidistribution of small points in case of algebraic dynamical systems over the function field of a curve. All these results assume $L$ ample, but we will see that in our case, the arguments hold more generally for big semiample line bundles. 

Yuan's variational principle for function fields will be proved in \S 5 and leads to the fundamental inequality between the height of $X$ and the heights $h_{\overline L}(P_m)$. It is based on Siu's theorem in the theory of big line bundles which we recall in \S 4. Note that Siu's inequality lies also at the root of Yuan's article \cite{Yu}. In fact, Yuan proves an arithmetic analogue in Arakelov geometry. In \S 6, we will first prove Theorem \ref{main theorem} from the fundamental inequality which is straightforward and which was used in all the articles mentioned above. Finally, we apply the methods to prove equidistribution of small subvarieties in algebraic dynamical systems.

%To have useful applications of Theorem \ref{main theorem}, a precise determination of the equidistribution measure $c_1(L,\metr)^{\wedge d}$ is indispensable. We restrict our attention now to an irreducible closed subvariety $X$ of the abelian variety $A$ over $K$ which may be eiter a function field or a number field. We assume that $L$ is an ample even line bundle on $A$ endowed with a canonical metric $\canmetr$.

To have useful applications of Theorem \ref{main theorem}, a precise determination of the equi\-distribution measure $c_1(L,\metr_v)^{\wedge d}$ is indispensable. We restrict our attention now to an irreducible $d$-dimensional closed subvariety $X$ of the abelian variety $A$ over $K$. Here, $K$ may be either a function field or a number field and we fix a non-archimedean place $v$ of $K$. We assume that $L$  is an even ample line bundle on $A$ endowed with a canonical metric $\canmetr$ (see Remark \ref{dynamical system for abelian variety} for the corresponding dynamical system). The height with respect to $(L,\canmetr)$ is the famous N\'eron--Tate height. 

In this situation, Theorem 6.7 of \cite{Gu6} gives a completely explicit description of the equidistribution measure $\mu:=c_1(L|_X,\metr_{{\rm can},v})^{\wedge d}$ in terms of convex geometry. More precisely, $\mu$ is supported in  finitely many rational simplices of dimension at most $d$ contained in $\Xan$ such that the restriction of $\mu$ to  each simplex is a positive multiple of the Lebesgue measure.

For applications, the following {\it tropical equidistribution theorem} is useful. We consider the Raynaud extension $1 \rightarrow T \rightarrow E \rightarrow B\rightarrow 0$ of $A_v^{\rm an}$. Here, $E$ is the uniformization of $A_v^{\rm an}$ such  that $A_v^{\rm an}=E/M$ for a lattice $M$ in $E$ and $B$ is an abelian variety over $\kdop$ of good reduction whose dimension is denoted by $b$. The split torus $T$ induces a map $\val:E \rightarrow \rdop^n$ with $\Lambda:=\val(M)$ a complete lattice in $\rdop^n$ (see \cite{Gu6}, \S 4). Passing to the quotient, we obtain a continuous map $\valbar:A_v^{\rm an} \rightarrow \rdop^n/\Lambda$. For a simplex $\Delta$ in $\rdop^n$, let $\delta_{\Deltabar}$ be the Dirac measure in $\Deltabar$, i.e. the push-forward of the relative Lebesgue measure on $\Delta$ to $\rtor$.

\begin{thm} \label{tropical equidistribution theorem}
Let $X$ be a geometrically integral $d$-dimensional closed subvariety of $A$. Then there is $e \in \{0, \dots, \min\{b,d\}\}$ such that the tropical variety $\valbar(\Xan)$ is the union of rational $(d-e)$-dimensional simplices $\overline{\Delta}_1, \dots ,\overline{\Delta}_M$ in $\rtor$ and there is (a possibly empty) list of rational simplices $\overline{\Delta}_{M+1}, \dots ,\overline{\Delta}_{M+N}$ contained in $\valbar(\Xan)$ with $\dim(\overline{\Delta}_j) \in \{d-b, \dots, d-e-1\}$ satisfying the following properties: 

For every ample even line bundle $L$ of $A$, there are  $r_1, \dots, r_{M+N} \in (0,\infty)$ such that for every generic net $(P_m)_{m \in I}$ in $X(\overline K)$ with 
$$\lim_m \hat{h}_L(P_m)=\frac{\hat{h}_L(X)}{(d+1)\deg_L(X)},$$ we have the following weak limit of regular probability measures on $\rtor$:
$$\frac{1}{|O(P_m)|} \sum_{P_m^\sigma \in O(P_m)} \delta_{\valbar(P_m^\sigma)} \stackrel{w}{\rightarrow} \sum_{j=1}^{M+N} r_j \delta_{\overline{\Delta}_j}$$
\end{thm}

This tropical equidistribution theorem follows immediately from Theorem \ref{main theorem} (resp. its number theoretic analogue in \cite{Yu}, Theorem 5.1) and from \cite{Gu6}, Theorem 1.1.

The skeleton $S(A_v^{\rm an})$ of $A$ is a canonical subset of $A_v^{\rm an}$ which is a proper deformation retraction of $A_v^{\rm an}$ and which is homeomorphic to $\rtor$ by the restriction of $\valbar$ (see \cite{Ber}, \S 6.5). Let $\mu$ be the unique regular Borel measure on $A_v^{\rm an}$ with $\supp(\mu)=S(A_v^{\rm an})$ such that $\valbar_*(\mu)$ is the Haar probability measure of $\rtor$.

\begin{thm} \label{equidistribution for X=A}
Let $K$ be either a function field or a number field and let $L$ be an ample even line bundle on the abelian variety $A$  over $K$. Let $(P_m)_{m \in I}$ be a generic net in $A(\overline K)$ such that $\lim_m \hat{h}_L(P_m)=0$. If $v$ is a non-archimedean place of $K$, then we have the following weak limit of regular probability measures on $A_v^{\rm an}$:
$$\frac{1}{|O(P_m)|} \sum_{P_m^\sigma \in O(P_m)} \delta_{P_m^\sigma} \stackrel{w}{\rightarrow} \mu.$$
\end{thm}

We will see in Theorem \ref{equidistribution of small subvarieties} and Remark \ref{dynamical system for abelian variety} that $\hat{h}_L(X)=0$. Then Theorem \ref{equidistribution for X=A} is a consequence of Theorem \ref{main theorem} (resp. \cite{Yu}, Theorem 5.1) and \cite{Gu6}, Corollary 7.3. Hence in the special case $X=A$, no lower dimensional simplices occur in Theorem \ref{tropical equidistribution theorem}. On the other hand, the quite natural Example 7.4 in \cite{Gu6} shows that lower dimensional simplices are possible in all dimensions $d-b, \dots, d-e-1$.

The author thanks the referee for his suggestions.

\vspace{3mm}
\centerline{\it Terminology}

In $A \subset B$, $A$ may be  equal to $B$. The complement of $A$ in $B$ is denoted by $B \setminus A$ \label{setminus} as we reserve $-$ for algebraic purposes. The zero is included in $\ndop$. We use $\overline K$ for an algebraic closure of a field $K$. 

A variety over $K$ is a separated reduced scheme of finite type. The group of cycles of pure dimension $t$ is denoted by $Z_t(X)$. If $L$ is a line bundle on a projective variety $X$, then $\deg_L(X)$ is the degree of $X$ with respect to $L$. If some positive tensor power of $L$ is generated by global sections, then $L$ is called {\it semiample}.

%The multiplicity of an irreducible component $Y$ of a scheme $S$ is denoted by $m(Y,S)$.

%Let $Y$ be a variety over a field. Following \cite{Ber4}, \S 2, we use the following canonical stratification of $Y$. We start with $Y^{(0)}:=Y$. For $r \in \ndop$, let $Y^{(r+1)} \subset Y^{(r)}$ be the complement of the set of normal points in $Y^{(r)}$. Since the set of normal points is open and dense, we get a chain of closed subsets:
%$$Y=Y^{(0)} \supsetneq Y^{(1)} \supsetneq Y^{(2)} \supsetneq \dots \supsetneq Y^{(s)} \supsetneq Y^{(s+1)} = \emptyset.$$
%The irreducible components of $Y^{(r)} \setminus Y^{(r+1)}$ are called the {\it strata} of $Y$. The set of strata is denoted by $\str(Y)$. It is ordered by $S \leq T$ if and only if $\overline{S} \subset \overline{T}$. A {\it strata subset} is a union of stratas. A {\it strata cycle} is a cycle whose components are strata subsets.

%For $\mb \in \zdop^n$, let $\xb^\mb:=x_1^{m_1} \cdots x_n^{m_n}$. The standard scalar product of $\ub,\ub' \in \rdop^n$ is denoted by  $\ub \cdot \ub':=u_1 u_1' + \dots + u_n u_n'$. For the notation used from convex geometry, we refer to \ref{convex geometry}. 

\section{Chambert-Loir's measures}

%Let $K=k(B)$ be the function field of an integral projective variety $B$ over the algebraically closed field $k$ such that $B$ is regular in codimension $1$. A place on $K$ is just a prime divisor and the valuation is the corresponding order and we denote the set of places by $M_B$. As usually, we fix an ample class $\mathbf c \in \Pic(B)$. If we count every prime divisor on $B$ with multiplicity $\deg_{\mathbf c}(Y)$, then we get a product formula on $K$ and hence a theory of heights which we now briefly describe in this section. 

In this section, $K_v$ is a field complete with respect to the discrete valuation $v$. For analytic considerations, we will work over the completion $\kdop$ of an algebraic closure of $K_v$.  This is an algebraically closed field with algebraically closed residue field (see \cite{BGR}, \S 3.4). The field $\kdop$ plays a similar role as $\cdop$ for archimedean places of number fields. The unique extension of $v$ to a valuation of $\kdop$ is also denoted by $v$. We fix a constant $c \in (0,1)$ and we use the absolute value $\absv:=c^{-v}$ to define {\it $v$-norms} $\metrv$ on a $\kdop$-vector space. We denote the valuation ring of $\kdop$ by $\kcirc$ and similarly we proceed for subfields of $\kdop$.

\begin{art} \rm \label{Berkovich analytic spaces}
In the following, we consider a projective scheme $X$ over $K_v$. We denote the {\it associated analytic space} by $\Xan$. Here, we use Berkovich's construction which behaves similar as the complex analytique analogue. Most algebraic properties hold also analytically and conversely, there is  a GAGA principle. For details, we refer to \cite{Ber}, 3.4.
\end{art}

\begin{art} \rm \label{v-metrics}
Let $L$ be a line bundle on $X$. A {\it $v$-metric} $\metrv$ on $L$ is a continuous family of $v$-norms on the fibres of $L_v^{\rm an}$. Two $v$-metrics $\metrv, \metrv'$ on $L$ induce a metric $\metrv'/\metrv$ on $O_X$ and evaluation at the constant section $1$ leads to a continuous nowhere vanishing function $g:=\metrv'/\metrv(1)$ on $\Xan$. The {\it distance of uniform convergence} is defined by
$$d_v(\metrv, \metrv'):=\sup_{x \in \Xan} |\log(g(x))|.$$
\end{art}

\begin{art} \rm \label{model and v-metrics}
A projective scheme $\Xcal$ over the valuation ring $K_v^\circ$ with generic fibre $X$ is called an {\it algebraic $K_v^\circ$-model} of $X$. If the line bundle $\Lcal$ on $\Xcal$ is an algebraic $K_v^\circ$-model of $L$, then we get a natural $v$-metric $\metr_\Lcal$ on $L$ by setting $\|s(x)\|_\Lcal=1$ for any local trivialization $s$ of $\Lcal$. We call $\metr_\Lcal$ the {\it algebraic $v$-metric} associated to $\Lcal$. 

More generally, we may consider models in the category of admissible formal schemes over  $\kcirc$  leading to {\it formal $v$-metrics} on $L$. A formal $v$-metric $\metr_\Lcal$ associated to the formal $\kcirc$-model $\Lcal$ is called {\it semipositive} if the reduction $\tilde{\Lcal}$ modulo $v$ is a numerically effective line bundle.
%Every line bundle has a {\it formal $v$-metric}. 
This formal point of view is suitable in the analytic context.
We refer to \cite{Gu2}, \S7, for more details. We will see in Proposition \ref{formal and algebraic} that every formal $v$-metric is algebraic over a finite base  change of $K_v$.

If $\metr_v$ is any formal $v$-metric on $O_X$, then $\log\|1\|_v$ is called a {\it formal function} on $\Xan$. The $\qdop$-subspace $\{\frac{1}{N} f\mid \text{$N \in \ndop\setminus \{0\}$,  $f$ formal function}\}$ is dense in $C(\Xan)$ (see \cite{Gu2}, Theorem 7.12). 

A $v$-metric $\metr_v$ on $L$ is called a {\it root of a formal $v$-metric} if some positive tensor power is a formal $v$-metric.
\end{art}

\begin{art} \rm \label{v-admissible metrics} \rm 
%For every $L$ and $X$ as above, there is a set $\ghpv(L)$ of continuous metrics on $L_v^{\rm an}$ with the following properties: These sets are closed  tensor product. If $\varphi:Y \rightarrow X$ is a morphism of projective schemes over $K_v$, then $\varphi^*\metrv \in \ghpv(\varphi^*L)$ if $\metrv \ghpv(L)$. The converse holds if $\varphi$ is surjective.  Every formal $v$-metric associated to a model $\Lcal$ with numerically effective reduction $\tilde{\Lcal}$ modulo $v$ is in $\ghpv(L)$. If $\metr_v^{\otimes m} \in \ghpv(L^{\otimes m})$ for some $m \in \ndop$, then $\metr_v \in \ghpv(L)$. The set $\ghpv(L)$ is closed under uniform convergence.

%We choose the smallest subsets $\ghpv(L)$ with these properties. Their elements are called {\it semipositive admissible $v$-metrics} on $L$.

A $v$-metric $\metr$ on $L$ is called a {\it semipositive admissible $v$-metric} if $\metr$ is the uniform limit of roots of semipositive formal $v$-metrics on $L$. A
$v$-metric $\metr_v$ of $L$ is called  {\it admissible} if there are  line bundles $L_1$, $L_2$ on $X$ with $\varphi^*(L) \cong L_1 \otimes L_2^{-1}$ such that $\metr_v = \metr_1/\metr_2$ for semipositive admissible   $v$-metrics $\metr_i$ of $L_i$.

We note that admissible $v$-metrics are closed under pull-back and tensor product. Every formal $v$-metric is admissible (see \cite{Gu3}, Proposition 10.4). Moreover, the canonical metrics on line bundles of an abelian variety over $K_v$ are admissible. This is proved by a variant of Tate's limit argument and is the very reason why we have allowed uniform limits and roots in the definition of semipositive admissible $v$-metrics (see \cite{Gu6}, Example 3.7). 

For details, we refer to \cite{Gu6}, \S 3. 
\end{art}

\begin{art} \rm \label{Chambert-Loir's measures for varieties}
In non-archimedean analysis, no good definition of Chern forms of $v$-metrized line bundles is known. However, Chambert-Loir \cite{Ch} has introduced a measure which is the analogue of top-dimensional wedge products of such Chern forms. Using a slight generalization of his construction and $d:=\dim(X)$, we get a regular Borel measure $c_1(\overline{L_1}) \wedge \cdots \wedge c_1(\overline{L_d})$ on $\Xan$ with respect to admissibly $v$-metrized line bundles $\overline{L_1},\dots, \overline{L_d}$, at least if $X$ is geometrically integral (see \cite{Gu6}, Proposition 3.8). Passing to a finite base extension and proceeding by linearity in the components, we deduce the following result:
\end{art}

\begin{prop} \label{Chambert-Loir's measures for cycles}
Let $L_1, \dots , L_t$ be line bundles on the projective scheme $X$ over $K_v$ endowed with admissible $v$-metrics $\metr_1, \dots, \metr_t$. For every cycle $Z \in Z_t(X)$ and every continuous function $g$ on $\Xan$, we have a real integral $\int_{\Zan} g \,c_1(\overline{L}_1) \wedge \dots \wedge c_1(\overline{L}_t)$ with the following properties:
\begin{itemize}
\item[(a)] The integral defines a bounded linear functional on $C(\Xan)$ and hence it is induced by a regular Borel measure on $\Xan$ with support in $\supp(Z)$. 
\item[(b)] The integral  is a multilinear symmetric function of $\overline{L_1}:=(L_1, \metr_1), \dots, \overline{L_t}:= (L_t,\metr_t)$ and linear in $Z$.
\item[(c)] If $\varphi:Y \rightarrow X$ is a morphism  projective schemes over ${K_v}$ and $W \in Z_t(Y)$, then 
$$\int_{\varphi_*(W)_v^{\rm an}} g \,c_1(\overline{L_1}) \wedge \dots \wedge c_1(\overline{L_t})
= \int_{W_v^{\rm an}} g \circ \varphi \, c_1(\varphi^*\overline{L_1}) \wedge \dots \wedge c_1(\varphi^*\overline{L_t}).$$
\item[(d)] If $\metr_1, \dots, \metr_t$ are semipositive, $Z \geq 0$ and $\mu:= c_1(\overline{L_1}) \wedge \dots \wedge c_1(\overline{L_t})$,   then
$$\left| \int_{\Zan} g \,\mu  \right| \leq \int_{\Zan} |g| \,\mu \leq |g|_{\rm sup} \deg_{L_1, \dots , L_t}(Z).$$
%\item[(e)] If $\Xcal$ is a formal $\kcirc$-model of $X$ with reduced special fibre $\tilde{\Xcal}$ and if $\metr_j$ is induced by a formal $\kcirc$-model $\Lcal_j$ of $L_j$ on $\Xcal$ for every $j=1, \dots, d$, then
%$$c_1(\overline{L_1}) \wedge \dots \wedge c_1(\overline{L_d})= 
%\sum_Y\deg_{\tilde{\Lcal}_1, \dots, \tilde{\Lcal}_d}(Y) \delta_{\xi_Y},$$
%where $Y$ ranges over the irreducible components of $\tilde{\Xcal}$ and $\delta_{\xi_Y}$ is the Dirac measure in the unique point $\xi_Y$ of $\Xan$ which reduces to the generic point of $Y$.
\item[(e)] We have $\int_{\Zan} c_1(\overline{L}_1) \wedge \dots \wedge c_1(\overline{L}_t)=\deg_{L_1,\dots,L_t}(Z)$. 
\item[(f)] For $j \in \{1, \dots, t\}$, let $\metr_{j,n}$ be a sequence of semipositive admissible $v$-metrics on $L_j$ converging uniformly to $\metr_j$. Then 
$$\lim_{n \to \infty} \int_{\Zan} g \,c_1(L_1,\metr_{1,n}) \wedge \dots \wedge c_1({L}_t, \metr_{t,n})
= \int_\Zan g \,c_1(\overline{L_1}) \wedge \dots \wedge c_1(\overline{L_t}).$$
%If we endow the set of positive regular Borel measures on $\Xan$ with the weak topology and if we fix the line bundles $L_1,\dots,L_d$ on $X$, then $c_1(\overline{L}_1) \wedge \dots \wedge c_1(\overline{L}_d)$ is continuous with respect to the vector $(\metr_1, \dots, \metr_d)$ of semipositive admissible metrics on $L_1,\dots ,L_d$.
%\item[(g)] If $X$ is a smooth variety and if there is a $j \in \{1, \dots,d\}$ such that $L_j$ is algebraically equivalent to $0$ endowed with a canonical metric, then $c_1(\overline{L}_1) \wedge \dots \wedge c_1(\overline{L}_d)=0$.
\end{itemize}
\end{prop}

%\proof We refer to \cite{Gu4}, \S 3, for the existence. Uniqueness is clear for formal metrics by (d). In general, it follows from a repeated application of the minimality of semipositive admissible metrics in \ref{v-admissible metrics} and will be skipped here. \qed

%\begin{rem} \label{Chambert-Loir for cycles}
%The restriction to geometrically integral varieties may be omitted by passing to a finite extension and then by linearity in the irreducible components. More generally, let $\overline{L_1}, \dots, \overline{L_t}$ be $v$-admissible metrized line bundles on a proper scheme $X$ over $K_v$. For a $t$-dimensional cycle $Z$ on $X$ and a continous function $f$ on $\Xan$, we define 
%$$\int_{Z_v^{\rm an}} f c_1(\overline{L}_1) \wedge \dots \wedge c_1(\overline{L}_t)$$
%by linearity in the geometrically irreducible components of $Z$.
%\end{rem}

\section{Heights}

Let $K=k(B)$ be the function field of an integral projective variety $B$ over the  field $k$ such that $B$ is regular in codimension $1$. A place $v$ on $K$ corresponds to a prime divisor $Y$ on $B$ and $v$ is the vanishing order along $Y$. By abuse of notation, we use $v$ also for the generic point of $Y$ and hence the set $M_B$ of places on $K$ may be viewed as a subset of $B$.
We fix an ample class $\mathbf c \in \Pic(B)$. If we count every $v \in M_B$ with multiplicity $\deg_{\mathbf c}(\overline v)$, then we get a product formula on $K$ for the discrete absolute values $\absv:=e^{-\ord_v(\cdot)}$ and hence a theory of heights which we now briefly describe in this section.

\begin{art} \label{admissible metrics} \rm 
Let $L$ be a line bundle on  a projective scheme $X$ over $K$. For $v \in M_B$, we apply the local theory from \S 2 to the line bundle $L_v:= L \otimes_K K_v$ over $X_v:= X \otimes_K K_v$. We define an {\it admissible $M_B$-metric} $\metr$ on $L$ as a family $\metr:=(\metr_v)_{v \in M_B}$ of admissible $v$-metrics $\metr_v$ on $L_v:=L \otimes_K K_v$ satisfying the following hypothesis:
There is an open dense subset $V$ of $B$, a projective scheme $\Xcal$ over $V$ with generic fibre $X$, $N \in \ndop \setminus \{0\}$ and a line bundle $\Lcal$ on $\Xcal$ such that $L^{\otimes N}= \Lcal|_X$ and such that $\metr_v^{\otimes N}$ is the algebraic $v$-metric associated to $\Lcal_v:= \Lcal \otimes_{K^\circ} K_v^\circ$ at all places $v \in M_B \cap V$.

The admissible $M_B$-metric $\metr$ on $L$ is called {\it semipositive} if $\metr_v$ is a semipositive admissible $v$-metric for all $v \in M_B$. Note that this definition of semipositivity is weaker than the one given in \cite{Gu5}, 3.2. 

It is clear that the pull-back and the tensor product of (semipositive) admissible $M_B$-metrics are again (semipositive) admissible $M_B$-metrics. We note that for any non-zero element $s$ in a fibre of $L$, there are only finitely many $v \in M_B$ with $\|s\|_v \neq 1$. 

The canonical metric of a line bundle on an abelian variety over $K$ is $M_B$-admissible. This is proved in \cite{Gu5}, 3.5, for ample symmetric line bundles. For odd line bundles, this follows similarly using \cite{Gu6}, Example 3.7, and the general case follows by linearity. The arguments will be generalized in \ref{dynamical systems} for arbitrary dynamical systems.
\end{art}

\begin{art} \rm \label{algebraic metrics}
An admissible $M_B$-metric $\metr$ on $L$ is called {\it algebraic} if $V=B$ and $N=1$ in the above definition. A {\it formal $M_B$-metric} is an admissible $M_B$-metric $\metr$ on $L$ such that $\metr_v$ is a formal $v$-metric for every $v \in M_B$. 
%A {\it model function} relative to $v \in M_B$ is given by $\log \|1\|_v$ for an algebraic $M_B$-metric on $O_X$.

%By \ref{model and v-metrics}, it is easy to see that every line bundle on a projective scheme $X$ over $K$ has a formal $M_B$-metric. If $X$ is projective, then 
Every line bundle $L$ is isomorphic to the difference of two very ample line bundles and hence we deduce that $L$ has an algebraic $M_B$-metric.
\end{art}

\begin{art} \label{base extension} \rm 
We are interested in the absolute height for points or more generally for cycles defined over an algebraic closure $\overline K$. Every finite subextension $K'/K$ is a function field $K'=k(B')$, where $B'$ is an irreducible normal projective variety with a finite surjective morphism $p:B'\rightarrow B$ leading to the finite extension $K'/K$ (see \cite{BG}, Lemma 1.4.10). 
%We consider the ample class $\mathbf c':=\phi^*(\mathbf c)$ on $B'$ and we count the prime divisors $Y'$ on $B'$ with multiplicity $\frac{1}{[K':K]} \deg_{\mathbf c'}(Y')$. 
On $K'$, we use the discrete absolute values extending those of $K$, i.e. for $w \in M_{K'}$ with ramification index $e(w/v)$ over $v \in M_K$, we consider the absolute value
$$|f|_w:=e^{-\ord_w(f)/e(w/v)}$$
on $K'$. We consider the ample class $\mathbf c':=p^*(\mathbf c)$ on $B'$ and we count every $w \in M_{K'}$ with multiplicity 
$$\mu(w):=\frac{e(w/v)\deg_{\mathbf c'}(\overline w)}{[K':K]} =\frac{[K'_w:K_v]\deg_{\mathbf c}(\overline v)}{[K':K]},$$ 
where $[K'_w:K_v]$ is the degree of the completions (see \cite{BG}, Example 1.4.13).

We show next that every admissible $M_B$-metric $\metr$ on $L$ induces a canonical $M_{B'}$-metric on $L':=L \otimes_K K'$ which we call the {\it base change} of $\metr$ to $L'$.  Every $w \in M_{B'}$ is lying over a unique $v \in M_{B}$. Using the identification $\kdop = {\mathbb K}_w$,  the $v$-metric on $L$ is a $w$-metric. The metric $\metr_w$ is induced by the base change of the algebraic model for $v$ in a dense open subset and hence we get a canonical $M_{B'}$-admissible metric on $L'$. If $\metr$ is semipositive, then it is clear that the base change of $\metr$ to $L'$ is also semipositive.
\end{art}

The following result clarifies the relation between formal and  algebraic metrics. The arguments are analogous to Yuan's Lemma 5.5. Since $\kcirc$ is not noetherian, additional care is needed here and we give the full proof.

\begin{prop} \label{formal and algebraic}
On a projective scheme $X$ over $K$, every formal metric is algebraic after a suitable finite base extension $K'/K$.
\end{prop}

\proof Let $\metr$ be an admissible $M_B$-metric on the line bundle $L$ of $X$. Since $L$ has an algebraic metric (see \ref{algebraic metrics}), we may assume that $L=O_X$. There is a dense open subset $V$ of $B$ such that $\metrv$ is the trivial metric on $O_X$ for every $v \in M_B \cap V$. For $v \not \in V$, $\metrv$ is the formal $v$-metric on $O_X$ induced by a formal $\kcirc$-model $\Lcal_v$ of $O_X$. The line bundle $\Lcal_v$ is living on a formal $\kcirc$-model $\Xcal_v$ of $\Xan$. We choose a projective $B$-model $\Ycal$ of $X$. By \cite{BL}, \S 4, we may assume that $\Xcal_v$ is the admissible formal blowing up of $\Ycal_v:=\Ycal \times_B \kcirc$ in a coherent ideal sheaf $\Jcal_v$ supported in the special fibre over $v$. By the formal GAGA principle (\cite{Ullr},  Theorem 6.8), $\Jcal_v$ is indeed defined algebraically and hence $\Xcal_v$ is the formal completion of a projective $\kcirc$-model $\Xcal_v^{\rm alg}$ along the special fibre. The formal GAGA principle again shows that $\Lcal_v$ is the formal completion of an algebraic $\kcirc$-model $\Lcal_v^{\rm alg}$ of $L \otimes_K \kdop$ on $\Xcal_v^{\rm alg}$. 

We choose a $K$-embedding $\sigma: \overline{K} \rightarrow \kdop$ and hence we get a place $u$ of $\overline K$. We note that $\kdop$ is the completion of $\overline{K}^\sigma:=\sigma(\overline{K})$ (\cite{BGR}, Proposition 3.4.1/3). Since $\Jcal_v$ is supported in the special fibre, $\Jcal_v$ contains a non-zero $\beta \in \overline{K}^\sigma$ with $|\beta|_v <1$. We easily deduce that $\Jcal_v$ may be generated by homogeneous polynomials with coefficients in $\overline{K}^\sigma \cap \kcirc$ and hence $\Xcal_v^{\rm alg}$ is defined over $\overline{K}^\sigma \cap \kcirc$. The line bundle $\Lcal_v^{\rm alg}$ is given by the Cartier divisor $D_v$ on $\Xcal_v^{\rm alg}$ associated to the meromorphic section $1$. Multiplying $D_v$ by a suitable non-zero element in $\overline{K}^\sigma$, we may assume that $D_v$ is an effective Cartier divisor on $\Xcal_v^{\rm alg}$. Then $O(-D_v)$ is  a coherent ideal sheaf on $\Xcal_v^{\rm alg}$ containing a sufficiently small $\beta' \in \overline{K}^\sigma \setminus \{0\}$. As above, we deduce that $O(-D_v)$ is a coherent ideal sheaf defined over $\overline{K}^\sigma \cap \kcirc$. Passing to a blowing up in this ideal sheaf supported in the special fibre over $v$, we conclude that $\metrv$ is induced by an algebraic model of $O_X$ defined over $\overline{K}^\sigma \cap \kcirc$. 

For a finite normal subextension ${K'}/K$ of $\overline{K}/K$, let $X' := X \otimes_K K'$, let $R_w$ be the valuation ring of the place $w:=u|_{K'}$ of ${K'}$  and let $\metr_w$ be the $w$-metric of $O_{X'}$ obtained from $\metr_v$ by base change.
Using $\sigma$-transfer of the above model to $\overline{K}$ and choosing ${K'}/K$ sufficiently large, we get an algebraic model $\Lcal_w'$ of $O_{X'}$ over  $R_w$ inducing $\metr_w$. Note that every place $w$ of ${K'}$ over $v$ is of this form for a suitable $\sigma$ and that ${K'}$ is independent of $\sigma$ by normality. Since there are only finitely many $v \in M_B \setminus V$, we may assume that ${K'}$ works for every such $v$.

Now ${K'}$ is the function field of an irreducible regular projective variety $B'$ over $k$ such that the extension ${K'}/K$ is induced by a finite morphism $p:B' \rightarrow B$. For a place $w$ in the dense open subset $ V':=p^{-1}(V)$ of $B'$, the metric $\metr_w$ is induced by the trivial bundle on $\Ycal':=\Ycal \times_B V'$. For $w \not \in V'$, the model $\Lcal_w'$ is defined and agrees with $\Ocal_{\Ycal'}$ over a dense open subset of $V'$. By glueing, we obtain an algebraic $B'$-model $\Lcal'$ of $O_{X'}$ inducing the base change of $\metr$ to $O_{X'}$. \qed

\begin{thm} \label{heights}
Let $K'=k(B')$ be a finite extension of $K$ as in \ref{base extension}. For admissibly $M_{B'}$-metrized line bundles $\overline{L_0}, \dots, \overline{L_t}$ on a projective scheme $X$ over $K'=k(B')$, there is a unique function $h_{\overline{L_0}, \dots, \overline{L_t}}:Z_t(X) \rightarrow \rdop$ with the following properties:
\begin{itemize} 
\item[(a)] $h_{\overline{L_0}, \dots, \overline{L_t}}(Z)$ is multilinear and symmetric in $\overline{L}_0 ,\dots ,\overline{L}_t$, and linear in $Z \in Z_t(X)$.
\item[(b)] If $\varphi:X' \rightarrow X$ is a morphism of projective schemes over $K'$, then we have functoriality
$$h_{\varphi^*\overline{L}_0, \dots, \varphi^*\overline{L}_t}=
 h_{\overline{L}_0 ,\dots ,\overline{L}_t}\circ \varphi_*.$$
\item[(c)] Let us replace the admissible $v$-metric $\metr_{0,v}$ at one place $v \in M_{B'}$ by the admissible $v$-metric $\metr_{0,v}'$ and let $\overline{L_0}'$ be the resulting $M_{B'}$-admissibly metrized line bundle. Then $\metr_{0,v}'/\metr_{0,v}$ is an admissible $v$-metric on $O_X$  defining a continuous function $g$ on $\Xan$ as in \ref{v-metrics} and we have
\begin{equation*}
\begin{split}
h_{\overline{L}_0',\overline{L}_1 ,\dots ,\overline{L}_t}(Z)
&- h_{\overline{L}_0 ,\dots ,\overline{L}_t}(Z) \\
&=
- \mu(v) \int_{Z_v^{\rm an}} \log(g) c_1(L_1,\metr_{1,v}) \wedge \dots \wedge c_1(L_t,\metr_{t,v}).
\end{split} 
\end{equation*}
\item[(d)] If the $M_{B'}$-metrics of $
\overline{L}_0 ,\dots ,\overline{L}_t$ are algebraic,   induced by the line bundles $\Lcal_0, \dots, \Lcal_t$ on the projective scheme $\pi:\Xcal \rightarrow B'$ with generic fibre $X$, then the height of any $t$-dimensional prime cycle $Z$ of $X$ with closure $\overline{Z}$ in $\Xcal$ is given as an intersection number on $\Xcal$ by
$$h_{\overline{L}_0 ,\dots ,\overline{L}_t}(Z)=\frac{1}{[K':K]}\deg_{\mathbf c'}\left(\pi_*\left(c_1(\Lcal_0) \dots c_1(\Lcal_t).\overline{Z} \right)\right).$$
\item[(e)] If $P \in X(K')$ and $\overline{L}$ is an admissibly $M_{B'}$-metrized line bundle on $X$, then
$$h_{\overline{L}}(P) = - \sum_{v \in M_{B'}} \mu(v) \log \| s(P) \|_v$$
for any non-zero $s(P)$ in the fibre of $L$ over $P$.
\end{itemize}
\end{thm}

\proof By multilinearity, we may assume that $L_0, \dots ,L_t$ are generated by global sections. In \S 11 of \cite{Gu3}, $h_{\overline{L}_0 ,\dots ,\overline{L}_t}(Z)$ was defined for any field with product formula. It is shown that (a),(b) and (e) hold. Property (c) is true by linearity in $\overline{L_0}$ and by definition of 
the right hand side as a local height of $Z$ (see \cite{Gu4}, 3.8). For algebraic metrics, the local heights in \cite{Gu3} are intersection numbers and (d) follows from the normalizations in \ref{base extension}. 

%By Chow's lemma and functoriality, we may reduce the general case to projective varieties. 
%Every line bundle is isomorphic to the ``difference'' of two line bundles generated by global sections. We may define $h_{\overline{L}_0 ,\dots ,\overline{L}_t}(Z)$ by multilinearity and (a)--(e) follows easily from the special case above.

%To prove uniqueness, we use that every line bundle has a formal metric (see \cite{Gu2},...). By (a) and (c), it is enough to prove uniqueness for formal metrics. 
%By functoriality (b) and Chow's lemma, it is enough to prove uniqueness for  projective varieties. Then 
Every line bundle has an algebraic $M_{B'}$-metric. By (a) and (c), we conclude that it is enough to consider algebraic $M_{B'}$-metrics and hence uniqueness follows from (d). \qed

\begin{Def} \rm \label{definition of height}
For $Z \in Z_t(X)$, we call $h_{\overline{L}_0 ,\dots ,\overline{L}_t}(Z)$ the {\it height} of $Z$ with respect to the $M_{B'}$-metrized line bundles ${\overline{L}_0 ,\dots ,\overline{L}_t}$. If $\overline{L}=\overline{L}_0= \dots =\overline{L}_t$, then the height of $Z$ is denoted by $h_{\overline L}(Z)$.
\end{Def}

\begin{rem} \rm \label{absolute height}
We note that $h_{\overline{L}_0 ,\dots ,\overline{L}_t}(Z)$ is an absolute height, i.e. it is invariant under base change and hence  well-defined on $Z_t(X \otimes_{K'}{\overline{K'}})$. To see this, let us consider the finite subextension $K''=k(B'')/K'$  of $\overline{K'}/K'$. Every $w \in M_{B''}$ is lying over a unique $v \in M_{B'}$. Then the $v$-metric on $L_j$ induces a canonical $w$-metric and hence we get an $M_{B''}$-admissible metric on $L_j \otimes_{K'} K''$. By (d) and the projection formula, the heights of $Z$ and $Z \otimes_{K'} K''$ agree with respect to algebraic $M_{B'}$-metrics. By uniqueness, they agree with respect to any $M_{B'}$-admissible metrics.

In particular, we get the height $h_{\overline{L}}(P)$ for every $P \in X(\overline{K'})$ with respect to an $M_{B'}$-metrized line bundle $L$ on $X$.
The distance of $M_{B'}$-metrics $\metr, \metr'$ on $L$ is measured by
$$d(\metr, \metr'):= \sum_{v \in M_{B'}} \mu(v)d_v(\metrv,\metr_v')$$ 
using the local distance $d_v$ of uniform convergence from \ref{v-metrics}. Since the metrics agree up to finitely many $v \in M_{B'}$, the sum is finite. By projection formula, the distance is absolute, i.e. invariant under base change of the metrics.
\end{rem}

Recall that Weil's theorem says that the height of points is determined  by the line bundle up to bounded functions. We have the following generalization for heights of cycles:

\begin{cor} \label{Weil's theorem}
For $j \in \{0, \dots, t\}$, let $\metr_j, \metr_j'$ be semipositive admissible $M_{B'}$-metrics on the line bundle $L_j$ of the projective scheme $X$ over $K'$. For every effective cycle $Z \in Z_t(X)$, we have
\begin{equation*}
\begin{split}
h_{(L_0,\metr_0'), \dots, (L_t,\metr_t')}(Z)&-h_{(L_0,\metr_0), \dots, (L_t,\metr_t)}(Z)\\
&\leq \sum_{j=0}^t d(\metr_j, \metr_j')\deg_{L_0, \dots, L_{j-1}, L_{j+1}, \dots , L_t}(Z)
\end{split}
\end{equation*}
and all occuring degrees are non-negative.
\end{cor}

\proof This follows from a $(t+1)$-fold application of Theorem \ref{heights}(d) and Proposition \ref{Chambert-Loir's measures for cycles}(c). \qed

\begin{ex} \rm \label{Neron-Tate}
Let $\overline{L}_0 ,\dots ,\overline{L}_t$ be canonically metrized line bundles on an abelian variety $A$ over $K$. We have seen in \ref{admissible metrics} that canonical metrics are admissible $M_B$-metrics. Since they are determined by the choice of a rigidification, they are unique up to positive rational  multiples. The product formula and Theorem \ref{heights} show that $h_{\overline{L}_0 ,\dots ,\overline{L}_t}(Z)$ does not depend on the choice of the canonical metrics. We call
$$\widehat{h}_{{L}_0 ,\dots ,{L}_t}(Z):=h_{\overline{L}_0 ,\dots ,\overline{L}_t}(Z)$$
the {\it N\'eron--Tate height} of $Z \in Z_t(A_{\overline{K}})$ with respect to $L_0, \dots, L_t$. In particular, we get a N\'eron--Tate height on $A(\overline{K})$ with respect to $L \in \Pic(A)$. We refer to \cite{Gu3} for the properties of $\widehat{h}_{{L}_0 ,\dots ,{L}_t}(Z)$ which can be easily deduced from Theorem \ref{heights}.
\end{ex}

\begin{art} \label{constant field extension} \rm  
All results obtained so far hold also for a number field $K$ using models over the ring of integers and arithmetic intersection theory. In the following, we consider special features of the function field $K=k(B)$. 

First, we deal with the extension of  the constant field $k$. Algebraic extensions are covered by \ref{base extension} and since the heights are absolute, we may assume in the following that $k$ is algebraically closed.

Let us consider a field extension $k'/k$ and the corresponding function field $K':=k'(B)$ of the variety $B':=B \otimes_k k'$. We will use always the ample class $\mathbf c'$ on $B'$ obtained from $\mathbf c$ by base change. We consider a line bundle $L$ on the projective scheme $X$ over $K$ endowed with an admissible $M_B$-metric $\metr=(\metr_v)_{v \in M_B}$. We claim that the line bundle $L':=L \otimes_K K'$ on $X':=X \otimes_K K'$ has a natural admissible $M_{B'}$-metric.

Indeed, there is an open dense subset $V$ of $B$ as in  \ref{admissible metrics}, where the metric $\metr$ is induced by the $N$-th root of a metric associated to a line bundle $\Lcal$ defined over $V$. By base change, we can use the model $\Lcal':=\Lcal \otimes_k k'$ on the open dense subset $V':=V \otimes_k k'$ to define $\metr_w':=\metr_{\Lcal',w}^{1/N}$ for all $w \in M_{B'} \cap V'$. Note that $M_{B'} \setminus V'$ is a finite set of prime divisors $v$ defined over $k$ and hence there is a unique extension of $\metr_v$ to a $v$-metric on $L'$. It is immediate from the definitions that $\metr'$ is $M_{B'}$-admissible.  

The same argument shows that $\metr'$ is semipositive if $\metr$ is semipositive. It follows from Theorem \ref{heights} and the compatibility of Chambert-Loir's measures with base change (\cite{Gu4}, Remark 3.10) that the height remains invariant under base change to $K'$.
\end{art}

\begin{art} \rm \label{def of generic curve}
In the remaining part of this section, we show how we can reduce heights over the function field $K=k(B)$ with $\delta:=\dim(B)\geq 2$ to the case of the function field of a curve. This is a classical process using a generic curve in $B$.

Again, we may assume that $k$ is algebraically closed. Replacing $\mathbf c$ by a tensor power, we may assume that $\mathbf c$ is very ample. We choose a basis $t_0, \dots ,t_N$ of global sections of a line bundle representing $\mathbf c$. Let  $\xi =(\xi_j^{(i)})$ with algebraically independent entries $\xi_j^{(i)}$, $i=1, \dots, \delta -1$, $j=0, \dots, N$, and let $\eta$ be the vector obtained from $\xi$ by omitting $\xi_0^{(i)}$, $i=1, \dots, \delta -1$ . Let $B'$ (resp. $B''$) be the base change of $B$ to $k':=k(\eta)$ (resp. $k'':=k(\xi)$). Then $B'$ and $B''$ are geometrically integral projective varieties which are geometrically regular in codimension $1$. The {\it generic curve} $B_{\mathbf c}$ is obtained by
$$B_{\mathbf c}:= \Div\left(\xi_0^{(1)}t_0+ \dots + \xi_N^{(1)}t_N\right) \dots \Div\left(\xi_0^{(\delta -1)}t_0+ \dots + \xi_N^{(\delta -1)}t_N\right).B''.$$
By \cite{Lan}, Section VII.6, $B_{\mathbf c}$ is a geometrically irreducible smooth projective curve over $k''$. By construction, we have 
\begin{equation} \label{function field identity}
k'(B')=k''(B_{\mathbf c}).
\end{equation}
Another way to see this is by considering the generic projection $\pi: B' \dashrightarrow \pdop_{k'}^{\delta -1}$, given by $$\pi[x_0: \dots :x_N]= [-x_0:\xi_1^{(1)} x_1 + \dots + \xi_N^{(1)}x_N: \dots : \xi_1^{(\delta-1)}x_1+ \dots + \xi_N^{(\delta-1)}x_N].$$
Then the fibre over the generic point $[1: \xi_0^{(1)}:\dots  :\xi_0^{(\delta-1)}]$ of $\pdop_{k'}^{\delta -1}$ is ${B_{\mathbf c}}$. Hence ${B_{\mathbf c}}$ is a dense subset of $B'$ with the same function field. Moreover, we get $M_{B_{\mathbf c}}=M_{B'} \cap {B_{\mathbf c}}$.
\end{art}

\begin{art} \rm \label{reduction to generic curve}
Now we apply \ref{constant field extension} and \ref{def of generic curve} to a line bundle $L$ on the projective scheme $X$ over $K$ endowed with an admissible $M_B$-metric $\metr$. We consider the base change to the function field $K':=k'(B')$, where $k'$ and $B'$ are as in \ref{def of generic curve}, leading to a line bundle $L':=L \otimes_K K'$ endowed with a natural metric $\metr'$ (see \ref{constant field extension}). By \eqref{function field identity}, $L'$ and $X'$ are defined over $k''({B_{\mathbf c}})$ and by restriction, we get a natural $M_{B_{\mathbf c}}$-metric $\metr_{\mathbf c}$ on $L'$. Since $\pi$ is a generic projection, the base change $w\in M_{B'}$ of $v \in M_B$ is contained in ${B_{\mathbf c}}$ and hence $w$ is a closed  point of ${B_{\mathbf c}}$ with valuation ring $\Ocal_{{B_{\mathbf c}},w}=\Ocal_{B',w}=\Ocal_{B,v} \otimes_k k'$ in $k''({B_{\mathbf c}})$. We claim that
\begin{equation} \label{generic height identity}
h_{({L}_0, \metr_0) ,\dots ,({L}_t,\metr_t)}(Z)= h_{({L}_0', \metr_{0,\mathbf c}) ,\dots ,({L}_t',\metr_{t, \mathbf c})}(Z')
\end{equation}
for every $Z \in Z_t(X)$ with $Z':=Z \otimes_K {K'}$ and line bundles $L_0, \dots , L_t$ with admissible $M_B$-metrics $\metr_0, \dots ,\metr_t$. 

\vspace{2mm}
\proof We deduce from \ref{constant field extension} that the height of $Z$ is equal to the height  of $Z'$ with respect to the function field of $B'$. We note that only places $w$ of $B'$ coming from places $v$ of $B$ contribute to the height of $Z'$ and hence we get \eqref{generic height identity}. More precisely, the height is given as a sum of local heights $\lambda(Z,v)$ with respect to suitable meromorphic sections $s_0, \dots , s_t$ (see \cite{Gu3}, \S 11). Over a dense open subset $V$ of $B$, $\lambda(Z,v)$ is given as an intersection number on an algebraic model. Using base change to the dense open subset $V':=V \otimes_k k'$ of $B'$, it is clear that $\lambda(Z,w)\neq 0$ only for places $w$ of $V'$ which are induced by places $v$ of $V$. Since the complement of $V'$ in $B'$ is also defined over $k$ and every irreducible component of $B' \setminus V'$ is the base change of an irreducible component of $B \setminus V$, we get the claim. \qed
\end{art}

\section{Big line bundles}

In this section, we recall some facts about big line bundles on a $t$-dimensional irreducible projective variety $\Xcal$ over any field $k$. For proofs, we refer to \cite{Laz}, \S 2.2. Note that Lazarsfeld states the results for projective varieties over an algebraically closed field of characteristic $0$, but the arguments hold in the more general setting.

\begin{Def} \label{big line bundles}
A line bundle $\Lcal$ on $\Xcal$ is called {\it big} if there is $C>0$ such that
$$h^0(\Lcal^{\otimes m}):= \dim H^0(\Xcal, \Lcal^{\otimes m}) \geq Cm^t$$
for all sufficiently large $m \in \ndop$.
\end{Def}

\begin{art} \label{pull-back and big} \rm 
Let $\varphi:\Xcal' \rightarrow \Xcal$ be a surjective morphism of irreducible projective varieties over $k$. We suppose that $\varphi$ is generically finite. If $\Lcal$ is big on $\Xcal$, then $\varphi^*(\Lcal)$  is big on $\Xcal'$. The converse holds for $\varphi$ birational. By passing to the normalization, we may always reduce to the case of normal varieties. 
\end{art}

The next result is called {\it Kodaira's lemma}.

\begin{lem} \label{Kodaira's lemma}
Let $D$ be a Cartier divisor on the irreducible projective variety $\Xcal$ over $k$. Then the following statements are equivalent:
\begin{itemize}
 \item[(a)] $O(D)$ is big.
 \item[(b)] For all ample divisors $H$ on $\Xcal$, there is a non-zero $m \in \ndop$ and an effective Cartier divisor $E$ on $\Xcal$ such that $mD$ is rationally equivalent to the Cartier divisor $H+E$.
 \item[(c)] There is an ample divisor $H$ on $\Xcal$, a non-zero $m \in \ndop$ and an effective Cartier divisor $E$ on $\Xcal$ such that $mD$ is numerically equivalent to $H+E$.
\end{itemize}
\end{lem}

\begin{art} \rm \label{exception set for big}
If $\Lcal$ is a big line bundle on the projective variety $\Xcal$, then there is a proper closed subset $\Vcal$ of $\Xcal$ such that for all irreducible closed subvarieties $\Ycal$ of $\Xcal$ not contained in $\Vcal$, the restriction $\Lcal|_\Ycal$ is big.
\end{art}

The main result of this section is {\it Siu's theorem}:

\begin{thm} \label{Siu's theorem} 
For nef Cartier divisors $D,E$ on $\Xcal$ and $m \in \ndop$, we have 
$$h^0(O(m(D-E))) \geq \frac{1}{t!}(D^t-tD^{t-1}\cdot E)m^t + o(m^{t}).$$
\end{thm}

Of course, this is most useful if the intersection number $D^t-tD^{t-1}\cdot E$ is positive. Then it follows that $O(D-E)$ is big. 
%In \cite{Laz}, Theorem 2.2.15, it is stated in this way (with $D,E$ nef instead of ample), but the proof gives Theorem \ref{Siu's theorem}.

\begin{thm} \label{nef and big}
Let $L$ be a numerically effective line bundle on $\Xcal$. Then $L$ is big if and only if $\deg_L(\Xcal)>0$. 
\end{thm}

In particular, the degree of $\Xcal$ with respect to a big semiample line bundle is positive.

\section{The fundamental inequality}

In this section, we prove the fundamental inequality which is a central tool to prove the equidistribution theorem in the next section. Through almost the whole section, we assume that $K$ is the function field of an irreducible regular projective curve $C$ over the algebraically closed field $k$.  Only at the end, we will generalize the fundamental inequality to arbitrary function fields using reduction to the generic curve.

We fix a big semiample line bundle $L$ on the $d$-dimensional irreducible projective variety $X$ over $K$.

\begin{Def} \rm \label{first successive minimum}
The {\it essential minimum} of $X$ with respect to the $M_C$-admis\-sibly metrized line bundle $\overline L$ is defined by 
$$e_1(X,\overline{L}):=\sup_Y \inf_{P \in X(\overline{K}) \setminus Y} h_{\overline{L}}(P),$$
where $Y$ ranges over all closed subsets of codimension $1$ in $X$. 
\end{Def}

\begin{art} \label{nef-metrics} \rm 
The semipositive admissible $M_C$-metric $\metr$ of $\overline L$ is called {\it nef} if $h_{\overline{L}}(P) \geq 0$ for all $P \in X(\overline K)$. If $\metr$ is an algebraic metric induced by a line bundle $\Lcal$ on a projective $C$-model $\Xcal$ of $X$, then this means that the restrictions of $\Lcal$ to all horizontal curves in $\Xcal$ have non-negative degrees. By semipositivity, the same holds for vertical curves in $\Xcal$ and hence $\metr_\Lcal$ is nef if and only if $\Lcal$ is nef.
\end{art}

\begin{lem} \label{twist of semipositive is nef}
Let $\metr_0$ be a semipositive admissible $M_C$-metric on $L$. Then there is an ample line bundle $\Mcal$ on $C$ such that $\metr_0 \otimes \pi^*\metr_\Mcal$ is a nef metric on $L$, where $\pi:X \rightarrow \Spec(K)$ is the morphism of structure.
\end{lem}

\proof Passing to a positive tensor power, we may assume that $L$ is generated by global sections and that there is an open dense affine subset $V$ of $C$ and a projective model $\pi:\Xcal \rightarrow V$ of $X$ with a line bundle $\Lcal$ on $\Xcal$ such that $\Lcal$ is a model of $L$ with $\metr_{0,v} = \metr_{\Lcal,v}$ for all $v \in  V$.  Let $s_1, \dots , s_n$ be a basis of $H^0(X,L)$. Using that $X=\Xcal \otimes_V K$ is obtained by flat base change, we get
$$H^0(\pi^{-1}(V),\Lcal) \otimes_{\Ocal(V)} K = H^0(X,L).$$
By passing to a smaller $V$, we may also assume that $s_1, \dots, s_n \in H^0(\pi^{-1}(V),\Lcal)$. For a closed point $v$ of $C$, $\|s_j\|_{0,v}$ is bounded. There is $R>0$ such that $\|s_j\|_{0,v} \leq R$ for all $v \in C \setminus V$ and $j=1, \dots ,n$. We consider the divisor
$$D:= \sum_{v \in C \setminus V} m\cdot [v]$$
for a multiplicity $m \geq \log(R)$. Then $s_D$ is a global section of $\Mcal:=\Ocal(D)$ with $m=-\log\|s_D\|_{\Mcal,v}$ for all $v \in C \setminus V$. For $P \in X(\overline K)$, there is $j \in \{1, \dots, n\}$ with $s_j(P) \neq 0$. There is an irreducible regular projective curve $C'$ and a finite morphism $p:C' \rightarrow C$ such that $P \in X(K')$ for the function field $K'=k(C')$. By Theorem \ref{heights}(e), we get
\begin{equation*}
\begin{split} 
h_{(L, \metr_0 \otimes \pi^*\metr_\Mcal})(P)&= - \sum_{w \in C'} \mu(w)\left(\log \|s_j(P)\|_{0,w} + \log \|s_D\|_{\Mcal, w} \right) \\
&\geq - \sum_{w \in C' \setminus p^{-1}(V)} \mu(w)   \left(\log \|s_j(P)\|_{0,w} - m \right) \geq 0
\end{split}
\end{equation*}
and hence $\metr_0 \otimes \pi^*\metr_\Mcal$ is nef. \qed 

\vspace{3mm} 
Using the distance $d_v$ of uniform convergence from \ref{v-metrics}, we have the following uniform version of the above result:

\begin{cor} \label{uniform twist}
Let $\metr_0$ be a semipositive admissible $M_C$-metric on $L$. Then there is  an ample line bundle $\Mcal$ on $C$ such that for every semipositive admissible $M_C$-metric $\metr$ on $L$, we get a semipositive admissible $M_C$-metric 
\begin{equation} \label{metric product}
\metr':=\metr \otimes \pi^*\metr_\Mcal \otimes \prod_{v \in M_C} \pi^* \metr_{O(v)}^{d_v(\metr_v,\metr_{0,v})} 
\end{equation}
on $L$ which is nef.
\end{cor}

\proof Since a positive tensor power of the metrics $\metr_0, \metr$ is given by models over an open dense subset of $C$, the metrics agree on a smaller open dense subset of $C$ and hence the product in \eqref{metric product} is finite. Using the notations from the proof of Lemma \ref{twist of semipositive is nef}, we have
\begin{equation*}
 \begin{split}
h_{(L, \metr')}(P) &= - \sum_{w \in C'} \mu(w)\left( \log\|s_j(P)\|_{w}+ \log\|s_D\|_{\Mcal, w} -d_w(\metr_w,\metr_{0,w}) \right)\\
& \geq - \sum_{w \in C'} \mu(w)\left( \log\|s_j(P)\|_{0,w}+ \log\|s_D\|_{\Mcal, w} \right) \geq 0
 \end{split}
\end{equation*}
and hence $(L, \metr')$ is nef. \qed 

\begin{art} \rm \label{uniformity wrt L} Let $\pi:\Xcal \rightarrow C$ be a model of $X$, i.e. $\pi$ is a projective flat morphism and $X$ is the generic fibre of $\Xcal$. Let $\Ncal$ be a line bundle on $\Xcal$ which is trivial on $X$. We consider models $\Lcal$ of $L$ on $\Xcal$ which are vertically nef, i.e. the restriction of $\Lcal$ to the fibre over $v$ is nef for all closed points $v$ of $C$. For $\ve \in \qdop$, we may view $\Lcal \otimes \Ncal^{\otimes \ve}$ as an element of $\Pic(\Xcal) \otimes_\zdop \qdop$ leading to well-defined heights and degrees using multilinearity. We will use $h_\Lcal:=h_{(L, \metr_\Lcal)}$ and $e_1(X,\Lcal):=e_1(X,(L,\metr_\Lcal))$. 

In the following lemma, we measure uniformity with respect to $\Lcal$ by $d(\Lcal):=d(\metr_\Lcal,\metr_0)$ using the distance to a fixed semipositive admissible $M_C$-metric $\metr_0$ on $L$ (see Remark \ref{absolute height}). 
\end{art}

\begin{lem} \label{Riemann-Roch inequality}
Under the hypothesis above and for $|\ve|\leq 1$, we have 
\begin{equation} \label{RR ineq}
\begin{split}
&h^0(\Xcal, (\Lcal \otimes \Ncal^{\otimes \ve})^{\otimes m})\\& \geq \frac{1}{(d+1)!} \left(\deg_{\Lcal \otimes \Ncal^{\otimes \ve}}(\Xcal) + (d(\Lcal)+1)O(\ve^2) \right) m^{d+1} + o_{\ve,\Lcal}(m^{d+1}) 
\end{split}
\end{equation}
for sufficiently large and sufficiently divisible $m \in \ndop$, where the implicit constant in $O(\ve^2)$ may depend on $L$ and $\Ncal$ but is independent of $\Lcal$, $m$ and $\ve$.
\end{lem}

\proof Replacing $\Ncal$ by $\Ncal^{-1}$ if necessary, we may assume that $\ve \geq 0$. There are semiample line bundles $\Lcal_1,\Lcal_2$ on $\Xcal$ with $\Ncal \cong \Lcal_1 \otimes \Lcal_2^{-1}$. As usual, we denote  the generic fibre of $\Lcal_i$ by $L_i$. First, we prove the lemma under the additional assumption that $\Lcal$ is nef. Then we have the decomposition 
\begin{equation}\label{L decomposition}
\Lcal \otimes \Ncal^{\otimes \ve} \cong \Lcal_+ \otimes (\Lcal_-)^{-1}
\end{equation}
for the numerically effective elements $\Lcal_+ := \Lcal \otimes \Lcal_1^{\otimes \ve}$ and $\Lcal_-:= \Lcal_2^{\otimes \ve}$ of $\Pic(\Xcal) \otimes \qdop$. For $m \in \ndop$ with $m \ve \in \zdop$, the $m$-th tensor power of \eqref{L decomposition} leads to a decomposition into numerically effective line bundles $\Lcal_+^{\otimes m}$ and $\Lcal_-^{\otimes m}$. We note that
$$c_1(\Lcal_+)^{d+1} - (d+1) c_1(\Lcal_+)^d \cdot c_1(\Lcal_-) = \deg_{\Lcal \otimes \Ncal^{\otimes \ve}}(\Xcal)+ O_\Lcal(\ve^2).$$
In fact, we may choose the implicit constant in $O_\Lcal(\ve^2)$  equal to
$$H_\Lcal:=3^{d+1}\max |h_{\underbrace{\Lcal,\dots,\Lcal}_{a},\underbrace{\Lcal_1, \dots, \Lcal_1}_{b}, \underbrace{\Lcal_2, \dots, \Lcal_2}_{c}}(X)|,
$$
where the maximum is taken over all natural numbers with $a+b+c=d+1$.  Since $\Lcal$, $\Lcal_1$ and $\Lcal_2$ are vertically nef, we may use Corollary \ref{Weil's theorem} to change from $\metr_\Lcal$ to $\metr_0$. Then we get
\begin{equation} \label{HL}
H_\Lcal \leq 3^{d+1} \left(H_0 + (d+1)d(\Lcal)D_L\right),
\end{equation}
where $ H_0$ is defined as $H_\Lcal$ with $(L,\metr_0)$ replacing $(L, \metr_\Lcal)$ and where 
$$D_L:=\max \deg_{\underbrace{L,\dots,L}_{a},\underbrace{L_1, \dots, L_1}_{b}, \underbrace{L_2, \dots, L_2}_{c}}(X),$$
with the maximum taken over all  $a+b+c=d$.
Hence we may apply Theorem \ref{Siu's theorem} to prove \eqref{RR ineq}.  

Now we prove the claim in general. For $v \in C$, we choose positive rational approximations $\delta_v \geq d_v(\metr_{\Lcal,v},\metr_{0,v})$. We may assume that 
\begin{equation} \label{good approximation}
\sum_{v \in C} \delta_v \leq d(\Lcal) + d(\Mcal),
\end{equation}
where $\Mcal$ is the ample line bundle on $C$ from Lemma \ref{twist of semipositive is nef} which is independent of $\Lcal$ and where $d(\Mcal):=d(\metr_\Mcal,\metr_{O_C})=\deg(\Mcal)$. 
By Corollary \ref{uniform twist},  $\Lcal':=\Lcal \otimes  \pi^*(\Mcal')$ is nef for
$$\Mcal':= \Mcal \otimes \left( \bigotimes_{v \in C} \Ocal_C(\delta_v[v]) \right).$$
For $m$ sufficiently large and sufficiently divisible, the first case shows
\begin{equation} \label{application of first case}
\begin{split}
&h^0(\Xcal, (\Lcal' \otimes \Ncal^{\otimes \ve})^{\otimes m})\\ &\geq \frac{1}{(d+1)!} \left(\deg_{\Lcal' \otimes \Ncal^{\otimes \ve}}(\Xcal) -H_{\Lcal'}\cdot\ve^2 \right) m^{d+1} + o_{\ve,\Lcal'}(m^{d+1}).
\end{split}
\end{equation}
Let $m_0 \in \ndop$ such that $(\Mcal')^{\otimes m_0}$ is very ample. Recursively, we choose a sequence $t_1, t_2, \dots $ of generic global sections of $(\Mcal')^{\otimes m_0}$. We conclude that the closed subschemes $E_i:=\pi^*(\Div(t_i))$ are disjoint on $\Xcal$. For sufficiently large and sufficiently divisible $m=m_0m_1$, we get an exact sequence
\begin{equation} \label{left exact sequence}
\begin{split}
0 \rightarrow H^0(\Xcal, (\Lcal \otimes \Ncal^{\otimes \ve})^{\otimes m}) \rightarrow &
H^0(\Xcal, (\Lcal' \otimes \Ncal^{\otimes \ve})^{\otimes m})\\
&\rightarrow \bigoplus_{i=1}^{m_1} H^0(E_i,(\Lcal' \otimes \Ncal^{\otimes \ve})^{\otimes m}),
 \end{split}\end{equation}
where the second map is induced by tensoring with $t_1 \otimes \dots \otimes t_{m_1}$ and the third map is induced by restriction. 
Since $\Ncal$ is trivial outside a finite set of fibres of $\Xcal$ over $C$, we conclude that $\Ncal|_{E_i}$ is trivial  for all  $i$ by the generic choice of $t_i$. Since $\pi^*(\Mcal')|_{E_i}$ is trivial for all $i$, we conclude that 
\begin{equation} \label{cohomology on the E_i}
H^0(E_i,(\Lcal' \otimes \Ncal^{\otimes \ve})^{\otimes m})=H^0(E_i,\Lcal^{\otimes m}).
\end{equation}
For $m$ sufficiently large and sufficiently divisible, we claim that 
\begin{equation} \label{independence of i} 
h^0(E_i,(\Lcal' \otimes \Ncal^{\otimes \ve})^{\otimes m}) \leq \frac{m_0}{d!}c_1(\Lcal)^d \cdot c_1(\pi^*\Mcal')m^d+o_\Lcal(m^{d})
\end{equation}
with an error term  independent of $i$. To prove this, we may assume $\Lcal'$ ample. Indeed, we may replace $\Lcal'$ by $\Lcal' \otimes \Hcal^{\otimes \lambda}$ for an ample line bundle $\Hcal$ on $\Xcal$ and rational $\lambda > 0$. The Nakai--Moishezon criterion shows that $\Lcal' \otimes \Hcal^{\otimes \lambda}$ is ample. Then continuity of the right hand side for  $\lambda \to 0$ proves \eqref{independence of i} in general. 

The ideal sheaf of $E_i$ is isomorphic to $\Ocal_\Xcal(-E_i)\cong \pi^*(\Mcal')^{-m_0}$ and hence we get the following long exact cohomology sequence:
\begin{equation*}
\begin{split}
0 \rightarrow H^0(\Xcal &, \pi^*(\Mcal')^{-m_0} \otimes {(\Lcal')}^{\otimes m})
\rightarrow H^0(\Xcal, {(\Lcal')}^{\otimes m}) \rightarrow H^0(E_i, {(\Lcal')}^{\otimes m})\\
\rightarrow& H^1(\Xcal, \pi^*(\Mcal')^{-m_0} \otimes {(\Lcal')}^{\otimes m})\rightarrow \cdots 
\end{split}
\end{equation*}
Since we may assume $\Lcal'$  ample, we have $H^1(\Xcal, \pi^*(\Mcal')^{-m_0} \otimes {(\Lcal')}^{\otimes m})=\{0\}$ for $m$ sufficiently large and we get independence of $h^0(E_i,(\Lcal')^{\otimes m})$ from $i$. By our above considerations in \eqref{cohomology on the E_i}, this dimension equals $h^0(E_i,(\Lcal' \otimes \Ncal^{\otimes \ve})^{\otimes m})$. 
Now we apply the Hilbert--Samuel formula in \eqref{cohomology on the E_i} to prove  \eqref{independence of i} and the independence of the error term from $i$. 

Using \eqref{independence of i} in \eqref{left exact sequence}, we get
\begin{equation} \label{first lower bound}
\begin{split} h^0(\Xcal, (\Lcal \otimes \Ncal^{\otimes \ve})^{\otimes m}) \geq &
h^0(\Xcal, (\Lcal' \otimes \Ncal^{\otimes \ve})^{\otimes m}) \\
&- \frac{1}{d!}c_1(\Lcal)^d \cdot c_1(\pi^*\Mcal')m^{d+1}+o_\Lcal(m^{d+1}).
\end{split}
\end{equation}
We may use the generic sections $t_i$ to compute the following intersection products
$$c_1(\pi^*\Mcal')^2=0, \quad c_1(\pi^*\Mcal'). c_1(\Ncal)=0$$
in $CH(\Xcal) \otimes_\zdop \qdop$ and hence we get
\begin{equation} \label{degree identity}
\deg_{\Lcal' \otimes \Ncal^{\otimes \ve}}(\Xcal) = \deg_{\Lcal \otimes \Ncal^{\otimes \ve}}(\Xcal) + (d+1)c_1(\Lcal)^d \cdot c_1(\pi^*\Mcal'). 
\end{equation}
For $v \in C$, we have
\begin{equation*}
\begin{split}
d_v(\metr_{\Lcal',v},\metr_{0,v})
&=d_v(\metr_{\Lcal,v},\metr_{0,v})+ d_v(\metr_{\Mcal',v},\metr_{\Ocal_C,v})\\
&=d_v(\metr_{\Lcal,v},\metr_{0,v})+ \delta_v+d_v(\metr_{\Mcal,v},\metr_{\Ocal_C,v})
\end{split}
\end{equation*}
and hence we get 
\begin{equation} \label{triangle inequality for metric}
d({\Lcal'}) \leq 2d(\Lcal)+2d(\Mcal)
\end{equation}
by \eqref{good approximation}. We deduce from \eqref{HL} and \eqref{triangle inequality for metric} that
\begin{equation} \label{HL'}
H_{\Lcal'}   \leq 3^{d+1} \left(H_0 + 2(d+1)\left(d(\Lcal)+d(\Mcal)\right)D_L\right).
\end{equation}
Using  \eqref{degree identity} and \eqref{HL'} in \eqref{application of first case}, we get a lower bound of $h^0(\Xcal, (\Lcal' \otimes \Ncal^{\otimes \ve})^{\otimes m})$ which we put in \eqref{first lower bound} to deduce \eqref{RR ineq}. \qed

\begin{rem} \rm \label{error term} 
The final argument shows that  the error term $(d(\Lcal)+1)O(\ve^2)$ in Lemma \ref{Riemann-Roch inequality} may be bounded  by $C_1\ve^2$ for
\begin{equation} \label{explicit error term formula}
C_1:=c_1(d)\left(H_0+\left(d(\Lcal)+d(\Mcal)\right)D_L\right).
\end{equation}
The constant $C_1=C_1(\Lcal,\Ncal,\metr_0)$ depends only on $(\Lcal,\Ncal,\metr_0)$ and the corresponding $\Mcal,\Lcal_1,\Lcal_2$ in the following way:
\begin{itemize}
 \item[(a)] For $m > 0$, we have $C_1(\Lcal^{\otimes m},\Ncal^{\otimes m},\metr_0^{\otimes m})=
m^{d+1}C_1(\Lcal,\Ncal,\metr_0)$. 
\item[(b)] $C_1$ is a universal positive continuous function in $(d,H_0,d(\Lcal),d(\Mcal),D_L)$ with $C_1=O_{d,D_L}(d(\Lcal)+d(\Mcal)+H_0)$.
\end{itemize}
\end{rem}

\begin{rem} \rm \label{absolute version}
We need an absolute version of Lemma \ref{Riemann-Roch inequality}. Let $K'$ be a finite extension of $K$. Then $K'$ is the function field of an irreducible regular projective curve $C'$ over $k$ such the extension $K'/K$ is induced by a finite morphism $p:C' \rightarrow C$. 

We still have the $C$-model $\Xcal$ of $X$ with the line bundle $\Ncal$ on $\Xcal$ which is trivial on the generic fibre $X$ and a fixed semipositive admissible $M_C$-metric $\metr_0$ of $L$. Now we consider a model $\Lcal'$ of $L':=L \otimes_K K'$, defined on a $C'$-model $\Xcal'$ of $X':=X \otimes_K K'$.  By Remark \ref{absolute height}, $d(\Lcal'):=d(\metr_{\Lcal'},\metr_0)$ is absolute, i.e. invariant under base change. 

There is always a $C'$-model $\Xcal''$ of $X'$ lying over $\Xcal$ and $\Xcal'$. Hence we may assume $\Xcal'=\Xcal''$ and that $\Ncal$ has a canonical pull-back $\Ncal'$ to $\Xcal'$. Then Lemma \ref{Riemann-Roch inequality}  holds with $\Xcal',\Lcal',\Ncal'$ replacing $\Xcal,\Lcal,\Ncal$ and with an absolute error term $(d(\Lcal')+1)O(\ve^2)$. Indeed, the invariants $d$, $H_0$, $d(\Lcal)$, $d(\Mcal)$ and $D_L$ are absolute and hence Remark \ref{error term} shows that the error term may be bounded by $C_1\ve^2$ with a constant $C_1$ independent of $K'$.
\end{rem}

\vspace{3mm}
In the following, we will use $C_2,C_3,\dots$ for constants with the properties (a) and (b) from Remark \ref{error term} and $c_2(d),c_3(d),\dots$ denote explicitly computable constants depending on $d$.
Now we are ready to prove the {\it fundamental inequality for algebraic metrics}. 

\begin{lem} \label{fundamental inequality for algebraic metrics}
We keep the assumptions and notations from \ref{uniformity wrt L}. There is a constant $c>0$ such that 
for all vertically nef $C$-models $\Lcal$ of $L$ on $\Xcal$ and all rational $\ve$ with $|\ve| \leq c$, we have $$\frac{h_{\Lcal \otimes \Ncal^{\otimes \ve}}(X)}{(d+1) \deg_L(X)} \leq e_1(X,\Lcal \otimes \Ncal^{\otimes \ve})+\left(d(\Lcal)+1\right)O(\ve^2),$$
where  the implicit constant in $O(\ve^2)$ is  independent of the choice of $\Lcal$ and $\ve$. More precisely, the constants and the inequality are absolute in the sense of Remark \ref{absolute version}.
\end{lem}

\proof We use the notation of the proof of Lemma \ref{Riemann-Roch inequality}. For $r \in \qdop$, we consider
$$\Lcal_r:= \Lcal \otimes \pi^*(\Mcal')^{\otimes r} \in \Pic(\Xcal) \otimes_\zdop \qdop.$$
For $m \in \ndop$ sufficiently large and sufficiently divisible, we would like to ensure the existence of a non-trivial global section $s$ of $(\Lcal_r \otimes \Ncal^{\otimes \ve})^{\otimes m}$ using Lemma \ref{Riemann-Roch inequality}. By Remark \ref{error term}, we have to choose $r$ such that 
\begin{equation} \label{choice of r}
%\begin{split}  
0  \stackrel{(!)}{<} h_{\Lcal_r \otimes \Ncal^{\otimes \ve}}(X) - C_1(\Lcal_r,\Ncal,\metr_0)\ve^2,
%&= h_{\Lcal \otimes \Ncal^{\otimes \ve}}(X) + r(d+1)\deg_{\Mcal'}(C) \deg_L(X)+ d(\Lcal_r)O(\ve^2),
%\end{split} 
\end{equation}
where we have used that the degree is equal to the height. The same computation as for \eqref{triangle inequality for metric} shows
$$d(\Lcal_r) \leq \left(|r|+1\right)d(\Lcal)+2|r|d(\Mcal)$$
and hence \eqref{explicit error term formula} yields 
\begin{equation} \label{C_1 for L_r}
C_1(\Lcal_r,\Ncal,\metr_0)\leq C_2|r|+C_1
\end{equation}
for $C_1:=C_1(\Lcal,\Ncal,\metr_0)$ and $C_2:=c_2(d)\left(d(\Lcal)+d(\Mcal)\right)D_L$. Similarly as in \eqref{degree identity} and by the projection formula, we have
\begin{equation} \label{height computation in choice of r}
h_{\Lcal_r \otimes \Ncal^{\otimes \ve}}(X)=h_{\Lcal \otimes \Ncal^{\otimes \ve}}(X) + r(d+1)\deg_{\Mcal'}(C)\deg_L(X).
\end{equation}
By \eqref{C_1 for L_r} and \eqref{height computation in choice of r}, the positivity assumption \eqref{choice of r} is satisfied if
\begin{equation} \label{necessary 1}
0<h_{\Lcal \otimes \Ncal^{\otimes \ve}}(X) +r\left((d+1)\deg_{\Mcal'}(C)\deg_L(X) \mp C_2\ve^2\right)-C_1\ve^2,
 \end{equation}
where the $``-"$ is used if and only if $r \geq 0$. Now we choose $\ve$  such that
\begin{equation} \label{sufficiently small}
C_2 \ve^2 \leq \frac{1}{2}(d+1)\deg_{\Mcal'}(C)\deg_L(X).
\end{equation}
By definition of $\Mcal'$,
$$R:=(d+1)\deg_{\Mcal'}(C)\deg_L(X)$$
is bounded below by
$$C_3:=(d+1)\left(d(\Lcal)+d(\Mcal)
\right) \deg_L(X)$$
and hence \eqref{sufficiently small}  holds for all
\begin{equation} \label{definition of c}
|\ve| \leq c:=\min\left\{\sqrt{\frac{C_3}{2C_2}},1\right\}=c_3(d)\cdot \sqrt{\frac{\deg_L(X)}{D_L}}.
\end{equation}
Moreover, we get
\begin{equation} \label{inverse}
\frac{1}{R} \leq \frac{1}{R-C_2\ve^2} \leq \frac{1}{R}\left(1+\frac{2C_2\ve^2}{R}\right) \leq 
\frac{1}{R} \left(1+\frac{\ve^2}{c^2}\right)
\end{equation}
and a similar estimate for $(R+C_2\ve^2)^{-1}$. The same arguments as in the first part of the proof of Lemma \ref{Riemann-Roch inequality} show that
\begin{equation} \label{absolute height bound}
|h_{\Lcal \otimes \Ncal^{\otimes \ve}}(X)| \leq |h_\Lcal(X)|+2^{d+1}H_\Lcal\ve
\leq c_4(d)\left(H_0+d(\Lcal)D_L\right)=:C_4.
\end{equation}
By \eqref{inverse} and \eqref{absolute height bound}, we conclude that \eqref{necessary 1} holds for
\begin{equation} \label{choice2 of r}
 r> \frac{-h_{\Lcal \otimes \Ncal^{\otimes \ve}}(X)+C_5\ve^2}{(d+1)\deg_{\Mcal'}(C) \deg_L(X)}, 
\end{equation}
where $C_5:=2C_1+C_4/c^2$. For such an $r$, we get a non-zero $s \in H^0(\Xcal,(\Lcal_r \otimes \Ncal^{\otimes \ve})^{\otimes m})$. Let $Y$ be the support of $\Div(s)$. For $P \in X(\overline K) \setminus Y$, there is a finite morphism $C' \rightarrow C$ of irreducible regular projective curves such that $P \in X(K')$ for the function field $K'=k(C')$. By Theorem \ref{heights}(e), we get
$$h_{\Lcal_r \otimes \Ncal^{\otimes \ve}}(P)=
-\frac{1}{m} \sum_{w \in C'} \mu(w) \log \|s(P)\|_{\Lcal_r \otimes \Ncal^{\otimes \ve},w} \geq 0$$ 
using that the metric of a global section of a model is bounded by $1$. We conclude that
\begin{equation} \label{first succ min non-negative}
e_1(X, \Lcal \otimes \Ncal^{\otimes \ve})+ r \deg(\Mcal') = e_1(X, \Lcal \otimes \Ncal^{\otimes \ve} \otimes \pi^*(\Mcal')^{\otimes r})\geq 0.
\end{equation}
If $r$ approaches the right hand side of \eqref{choice2 of r} and if we put this into \eqref{first succ min non-negative}, then we get 
\begin{equation} \label{precise fundamental inequality for algebraic metrics}
\frac{h_{\Lcal \otimes \Ncal^{\otimes \ve}}(X)}{(d+1) \deg_L(X)} \leq e_1(X,\Lcal \otimes \Ncal^{\otimes \ve})+S\ve^2
\end{equation}
for $S:= C_5/\{(d+1)\deg_L(X)\}$. Since the constants $C_1,C_2,\dots$ and $c$ are absolute, we deduce as in Remark \ref{absolute version} that the fundamental inequality is absolute. \qed

\vspace{3mm}

Recall that $L$ is a big semiample line bundle on $X$. Again, $\metr_0$ denotes a fixed semipositive admissible $M_C$-metric on $L$ and $\metr_f$ is a formal $M_C$-metric on $O_X$. 

%We fix a place $v \in M_C$ and we consider a formal function $f$ relative to $v$. Then $f$ is a continuous function on $X_v^{\rm an}$ induced by a formal $v$-metric $\metr_{f,v}$ of $O_X$. Let $\overline L(\ve f)$ be the line bundle $L$ endowed with the following admissible $M_B$-metric. If $w \in M_B$, $w \neq v$, then we use $\metr_w$ from $\overline L$. For $w=v$, we use the $v$-metric $\metr_v \otimes \metr_{f,v}^{\otimes \ve}$.
%\end{rem}

The following variational form of the {\it fundamental inequality} was proved by Yuan in the number field case (with $L$ ample, see \cite{Yu}, 5.2 and 5.3). Note that the absolute version of the error term is new here.

\begin{prop} \label{fundamental inequality}
There is an absolute constant $c>0$ with the following property: For every $\ve \in (-c,c)$ and every semipositive admissible $M_C$-metric $\metr$ of $L$, we have
\begin{equation} \label{fundamental inequality displayed}
\frac{h_{(L,\metr \otimes \metr_f^{\otimes \ve})}(X)}{(d+1) \deg_L(X)} \leq e_1(X,(L,\metr \otimes \metr_f^{\otimes \ve}))+(d(\metr, \metr_0)+1)O(\ve^2),
\end{equation}
where the implicit constant in $O(\ve^2)$ may depend on $L$ and $\metr_f$ but is independent of $\metr$ and $\ve$. Moreover, the implicit constant is absolute, i.e. it holds also for all semipositive admissible $M_{C'}$-metrics on $L \otimes_K K'$ independently of the finite extension $K'=k(C')$ of $K$.
\end{prop}

\proof If we use $d(\overline L):=d(\metr,\metr_0)$ instead of $d(\Lcal)$, then property (b) in Remark \ref{error term} shows that the constants $C_1,C_2, \dots$ and $S$ from our previous results make sense also for $\overline L=(L,\metr)$. We will prove that the fundamental inequality 
\begin{equation} \label{precise fundamental inequality for semipositive metrics}
\frac{h_{(L,\metr \otimes \metr_f^{\otimes \ve})}(X)}{(d+1) \deg_L(X)} \leq e_1(X,(L,\metr \otimes \metr_f^{\otimes \ve}))+S \ve^2 
\end{equation}
holds for $\overline L$ and $\ve \in (-c,c)$. This  implies the absolute nature of the fundamental inequality. Passing to a finite base extension, we may assume that $\metr_f$ is algebraic (Proposition \ref{formal and algebraic}), i.e. there is a projective flat scheme $\Xcal$ over $C$ with generic fibre $X$ and a line bundle $\Ncal$ on $\Xcal$ which is trivial on $X$ such that $\metr_{f}=\metr_{\Ncal}$.  

First, we prove the claim if $\metr$ is a formal metric and $\ve \in \qdop$. By Proposition \ref{formal and algebraic}, there is a finite extension $K'=k(C')/K$ such that the base change $\metr'$ of $\metr$ to $K'$ is an algebraic metric and the claim follows from \eqref{precise fundamental inequality for algebraic metrics}. 
%The absolute version of Lemma \ref{fundamental inequality for algebraic metrics} given in Remark \ref{absolute version} proves now \eqref{fundamental inequality displayed} with an absolute error term.

Next, we prove that \eqref{precise fundamental inequality for semipositive metrics} holds for roots of formal metrics and $\ve \in \qdop$. Indeed, the definition of $c$ in \eqref{definition of c} shows that $c$ depends only on $(L,\Ncal)$ and this homogeneously of degree $0$. Using property (a) of the constant $C_5$, we easily deduce that $S$ from \eqref{precise fundamental inequality for algebraic metrics} is a homogeneous function of $(\Lcal,\Ncal,\metr_0)$ of degree $1$.  This allows us to deduce \eqref{precise fundamental inequality for semipositive metrics} for roots of formal metrics from the case of formal metrics.

In general, $\metr$ is the limit of roots of semipositive formal $M_C$-metrics $\metr_n$ with respect to the distance  of metrics on $L$. We choose rational approximations $\ve_n$ of $\ve$. By the special case above, inequality \eqref{precise fundamental inequality for semipositive metrics} holds for $(L,\metr_n)$, $\ve_n$ and $S_n$ instead of $(L,\metr)$, $\ve$ and $S$. By Theorem \ref{heights}(e), we deduce easily
\begin{equation} \label{uniformity of successive minima}
\lim_{n \to \infty} e_1(X,(L,\metr_n \otimes \metr_f^{\otimes \ve_n} ))=e_1(X,(L,\metr \otimes \metr_f^{\otimes \ve})).
\end{equation}
Since $\Xcal$ is projective, we may write $\Ncal$ as the ''difference'' of two ample line bundles on $\Xcal$. Using multilinearity of the height and Corollary \ref{Weil's theorem}, we get
\begin{equation} \label{uniformity of the height of X}
\lim_{n \to \infty} h_{(L,\metr_n \otimes \metr_f^{\otimes \ve_n}) }(X) = h_{(L,\metr \otimes \metr_f^{\otimes \ve})}(X).
\end{equation}
Since $c$ is independent of the metrics and $S_n \to S$ (by property (b) in Remark \ref{error term}), we deduce \eqref{precise fundamental inequality for semipositive metrics} for $(L,\metr)$ immediately from \eqref{uniformity of successive minima} and \eqref{uniformity of the height of X}.  \qed

\begin{prop} \label{fundamental inequality for arbitrary function fields}
The fundamental inequality in the form given in Proposition \ref{fundamental inequality} holds for arbitrary function fields. 
\end{prop}

\proof  The heights and the essential minimum are invariant under algebraic base extension of $K$  and hence we may assume that $K$ is the function field of an irreducible projective variety $B$ over the algebraically closed field $k$ such that $B$ is regular in codimension $1$. Clearly, we may assume that $\mathbf c$ is very ample. Let $C:=B_{\mathbf c}$ be the generic curve over $k''$ constructed in \ref{def of generic curve}. We have seen in \ref{reduction to generic curve} that base change to the function field $K':=k''(B_{\mathbf c})$ induces a canonical semipositive admissible $M_C$-metric $\metr_{\mathbf c}$ on $L':=L \otimes_K K'$ such that heights are invariant under base change to $K'$. Since the fundamental inequality holds for $K'$ by Proposition \ref{fundamental inequality} and since the distance between the metrics is invariant under base change, it is enough to prove that the essential minimum decreases under base change to $K'$.

Let $Y'$ be a proper closed subset of $X'=X \otimes_K K'$. Using a suitable specialization, it is easy to construct a proper closed subset $Y$ of $X$ such that $P \not \in Y(\overline K)$ implies $P \not \in Y'$. We get
$$\inf_{P \in X(\overline K) \setminus Y} h_{(L,\metr)}(P)
\geq \inf_{P \in X(\overline K) \setminus Y'} h_{(L,\metr)}(P)
\geq \inf_{P \in X'({\overline{ K'}}) \setminus Y'} h_{(L,\metr_{\mathbf c})}(P)$$
proving immediately that the essential minimum decreases under base change to $K'$. \qed

\section{Equidistribution theorems}

In this section, $K$ is always a function field. First, we will deduce Theorem \ref{main theorem} from the fundamental inequality. Then we apply our methods to dynamical systems. In this case, we will generalize equidistribution to small generic nets of subvarieties. As examples, we will consider abelian varieties and multiplicative tori.

\vspace{3mm}
\noindent {\bf Proof of Theorem \ref{main theorem}:\:}
We have seen in \ref{model and v-metrics} that formal functions generate a dense $\qdop$-subspace of $C(\Xan)$, hence it is enough to check 
\begin{equation} \label{equidistribution for formal functions}
\lim_m \frac{1}{|O(P_m)|} \sum_{P_m^\sigma \in O(P_m)} f(P_m^\sigma)
= \frac{1}{\deg_L(X)} \int_{\Xan} f c_1(L,\metr_v)^{\wedge d}
\end{equation}
for a formal function $f$ on $\Xan$. Then $f=-\log\|1\|_{f,v}$ for a formal $M_B$-metric $\metr_f$ of $O_X$ which is trivial for all places different from $v$. Since the net $(P_m)_{m \in I}$ is generic in $X$, the fundamental inequality (Proposition \ref{fundamental inequality for arbitrary function fields}) yields
\begin{equation} \label{application of the fundemental inequality}
\frac{1}{(d+1)\deg_L(X)} h_{(L, \metr \otimes \metr_f^{\otimes \ve})}(X) \leq
\liminf_m h_{(L, \metr \otimes \metr_f^{\otimes \ve})}(P_m) + O(\ve^2), 
\end{equation}
where the implicit constant is independent of $\ve \in (-c,c)$. By Theorem \ref{heights}, we easily deduce
$$h_{(L, \metr \otimes \metr_f^{\otimes \ve})}(X) =h_{\overline L}(X)+ \ve(d+1)\mu(v)\int_\Xan f c_1(L,\metr_v)^{\wedge d} + O(\ve^2)$$
and 
$$h_{(L, \metr \otimes \metr_f^{\otimes \ve})}(P_m)= h_{\overline L}(P_m) + 
\ve\frac{\mu(v)}{|O(P_m)|} \sum_{P_m^\sigma \in O(P_m)} f(P_m^\sigma) + O(\ve^2).$$
If we insert both identities in \eqref{application of the fundemental inequality}, then $\lim_m h_{\overline L}(P_m)=\frac{h_{\overline L}(X)}{(d+1)\deg_L(X)}$ leads to
$$ \frac{\ve}{\deg_L(X)} \int_{\Xan} f c_1(L,\metr_v)^{\wedge d}
\leq  \liminf_m \frac{\ve}{|O(P_m)|} \sum_{P_m^\sigma \in O(P_m)} f(P_m^\sigma) + O(\ve^2).$$
Choosing $\ve >0$ and letting $\ve \to 0$, we get $``\geq"$ in \eqref{equidistribution for formal functions} with  $``\lim"$ replaced by $``\liminf"$. Using the same argument for $-f$ instead of $f$, we get also $``\leq"$ for the $``\limsup"$ and hence the claim. \qed

\begin{art} \label{dynamical systems} \rm 
We consider a dynamical system $(W,\Phi,L)$, where $L$ is an ample line bundle on the irreducible projective variety $W$ over the function field $K=k(B)$ and $\Phi:W \rightarrow W$ is a morphism with $$\Phi^*(L)\cong L^{\otimes q}$$ for some $q \in \ndop$ with $q >1$. This is an isometry with respect to a unique admissible $M_B$-metric $\canmetr$ of $L$. In fact, the canonical metric is given as the uniform limit of recursively defined metrics 
$$\metr_j :=\left(\Phi^*\metr_{j-1} \right)^{1/q}$$
on $L$, where $\metr_0$ is any admissible $M_B$-metric on $L$ (see \cite{BG}, Theorem 9.5.4 and its proof). We may choose $\metr_0$ as a root of an algebraic metric associated to an ample model and hence $\canmetr$ is semipositive, even uniform limit of roots of ample metrics.  We call 
$$\hat{h}_L(Z):= h_{(L,\canmetr)}(Z)$$
the {\it canoncial height} of the cycle $Z$ of $W$. Since $\canmetr$ is uniform limit of roots of ample metrics, the canonical height is non-negative for effective cycles. By functoriality in Theorem \ref{heights}, we deduce easily $\hat{h}_L(W)=0$. 
\end{art}

\begin{art} \rm \label{setting for equidistribution of small subvarieties}
In the number field case, Yuan (\cite{Yu}, Theorem 5.6) proved  a general {\it equidistribution theorem of small subvarieties} motivated by the case of abelian varieties due to Autissier \cite{Au} and Baker--Ih \cite{BI}. Next, we will give a similar equidistribution theorem for the dynamical system $(W,\Phi,L)$ over the function field $K=k(B)$. A similar result was proved independently by Faber (\cite{Fa}, Theorem 4.1) in the case of function fields of curves.

Let $X$ be an irreducible $d$-dimensional closed subvariety of $W$. We consider a {\it generic} net $(Y_m)_{m \in I}$ of irreducible $t$-dimensional closed subvarieties of $X \otimes_K \overline K$, i.e. for every proper closed subset $Y$ of $X$, there is $m_0 \in I$ such that $Y_m$ is not contained in $Y$ for all $m \geq m_0$. We  assume that the net is {\it small} in the sense  of
\begin{equation} \label{small for subvarieties}
\lim_m \frac{1}{(t+1)\deg_L(Y_m)} \hat{h}_L(Y_m)=0.
\end{equation}
Let us denote the $\Gal(\overline K/K)$-orbit of $Y_m$ by $O(Y_m)$. We fix a place $v \in M_B$ and an embedding $\overline K \hookrightarrow \kdop$ over $K$ to identify $Y_m \otimes_K \overline K$ with a subvariety of $\Xan$. Then we may view $\deg_L(Y_m)^{-1}c_1(L|_{Y_m},\metr_{{\rm can},v})^{\wedge t}$ as  a regular probability measure on $\Xan$ with support in $(Y_m)_v^{\rm an}$. 
\end{art}

\begin{thm} \label{equidistribution of small subvarieties}
Under the hypothesis of \ref{equidistribution of small subvarieties}, we have $\hat{h}_L(X)=0$ and the following weak limit of regular probability measures on $\Xan$:
\begin{equation*} \label{weak limit for small subvarieties}
\frac{1}{|O(Y_m)|\deg_L(Y_m)} \sum_{Y_m^\sigma \in O(Y_m)} c_1(L|_{Y_m^\sigma},\metr_{{\rm can},v})^{\wedge t} 
\stackrel{w}{\rightarrow} \frac{1}{\deg_L(X)} c_1(L|_X,\metr_{{\rm can},v})^{\wedge d}.
\end{equation*}
\end{thm}

\proof 
%Since $\canmetr$ is the uniform limit of roots of ample metrics, Zhang's theory of successive minima applies. By \cite{Gu5}, \S 4, this holds also in the function field case and \cite{Gu}, Proposition 4.3 shows
%$$\deg_L(X)e_1(Y_m,(L,\canmetr)) \leq \hat{h}(Y_m).$$
%Using \eqref{small for subvarieties} and that the net is generic in $X \otimes_K \overline K$, we conclude $e_1(X,(L,\canmetr))=0$. The fundamental inequality (Proposition \ref{fundamental inequality for arbitrary function fields}) yields now $\hat(X)=0$.
We need some preliminary considerations using the methods of \S 5. We first assume that $K$ is the function field of an irreducible projective curve $C$ over an algebraically closed field $k$ and that $\mathbf c$ has degree $1$. We also assume that $\metr$ is an algebraic $M_C$-metric on $L$  induced by an ample model $\Lcal$ on a flat projective scheme $\pi:\Xcal \rightarrow C$ with generic fibre $X$. As in \ref{uniformity wrt L}, we consider a line bundle $\Ncal$ on $\Xcal$ with $\Ncal|_X=O_X$ and $\ve \in \qdop$. We use the notation of \ref{Riemann-Roch inequality}--\ref{fundamental inequality for algebraic metrics}. For $\ve \in (-c,c)$, there is a non-zero global section $s$ of $(\Lcal_r \otimes \Ncal^{\otimes \ve})^{\otimes m}$ if 
\begin{equation} \label{condition for r}
r> \frac{-h_{\Lcal \otimes \Ncal^{\otimes \ve}}(X)+C_5\ve^2}{(d+1)\deg_{\Mcal'}(C) \deg_L(X)}
%+\frac{d(\Lcal)}{\deg_{\Mcal'}(C)} O(\ve^2)
\end{equation}
and if $m$ is sufficiently large and sufficiently divisible (see \eqref{choice2 of r}). Let $Y$ be a generic $t$-dimensional closed subvariety of $X \otimes_K \overline K$ and $Z:=\Div(s|_Y) \in Z_{t-1}(X)$. Let $K'$ be a finite extension of $K$ such that $Y$ is defined over $K'$ and let $X',L'$ be the base change of $X,L$ to $K'$. By the induction formula for heights of subvarieties (see \cite{Gu3}, Remark 9.5) or by an easy direct calculation of intersection numbers, we get
\begin{equation} \label{induction formula applied}
\begin{split} 
&h_{\Lcal, \dots, \Lcal,(\Lcal_r \otimes \Ncal^{\otimes \ve})^{\otimes m}}(Y) \\
&=h_\Lcal(Z) - \sum_{w \in M_{K'}}\mu(w) \int_{{Y_w^{\rm an}}} \log\|s\|_w' c_1(L'|_Y,\metr_{\Lcal,w})^{\wedge t},
\end{split} 
\end{equation}
where $\metr_w'$ is induced by the model  $(\Lcal_r \otimes \Ncal^{\otimes \ve})^{\otimes m}$. Since $\|s\|_w'\leq 1$, we get
\begin{equation*}
h_\Lcal(Y) + r \deg_{\Mcal'}(C)\deg_L(Y)+ \ve h_{\Lcal, \dots ,\Lcal,\Ncal}(Y) \geq \frac{1}{m} h_\Lcal(Z).
\end{equation*}
Using  $h_\Lcal(Z) \geq 0$ for the effective cycle $Z$, we deduce
$$\frac{1}{\deg_L(Y)} h_\Lcal(Y) \geq -r \deg_{\Mcal'}(C) - \ve \frac{h_{\Lcal, \dots ,\Lcal,\Ncal}(Y)}{\deg_L(Y)}.$$
Now the right hand side is also independent of $s$. If $r$ approaches the right hand side of \eqref{condition for r}, then we get
$$\frac{1}{\deg_L(Y)} h_\Lcal(Y) \geq \frac{h_{\Lcal \otimes \Ncal^{\otimes \ve}}(X)-C_5\ve^2}{(d+1)\deg_L(X)} - \ve \frac{h_{\Lcal, \dots ,\Lcal,\Ncal}(Y)}{\deg_L(Y)}.$$
If we use Corollary \ref{Weil's theorem} for $h_{\Lcal \otimes \Ncal^{\otimes \ve}}(X)$ in the same way as in the proof of Lemma \ref{Riemann-Roch inequality}, then the right hand side has the lower bound
\begin{equation*} \label{algebraic fundamental inequality for Y}
\frac{h_{\Lcal}(X)}{(d+1)\deg_L(X)}
+\ve \left( \frac{h_{\Lcal, \cdots ,\Lcal, \Ncal}(X)}{\deg_L(X)}-\frac{h_{\Lcal, \dots ,\Lcal,\Ncal}(Y)}{\deg_L(Y)}\right) -\frac{C_6\ve^2}{\deg_L(X)}.
\end{equation*}
Similarly as for the fundamental inequality, we can extend the resulting inequality to arbitrary semipositive admissible $M_C$-metrics on $L$ and to any formal metric $\metr_f$ on $O_X$ to get
\begin{equation} \label{semipositive fundamental inequality for Y}
\begin{split}
\frac{1}{\deg_L(Y)} 
h_{\overline L}(Y) 
\geq &\frac{h_{\overline L}(X)}{(d+1)\deg_L(X)} \\
&+\ve \left( \frac{h_{{\overline L}, \cdots ,{\overline L}, \overline O_X^f}(X)}{\deg_L(X)}-\frac{h_{{\overline L}, \dots ,{\overline L},\overline O_X^f}(Y)}{\deg_L(Y)}\right) -\frac{C_6\ve^2}{\deg_L(X)},
\end{split}
\end{equation}
where $\overline O_X^f:=(O_X,\metr_f)$. Again, the constant $C_6$ is absolute. Next, we note that \eqref{semipositive fundamental inequality for Y} holds for arbitrary function fields $K=k(B)$. Indeed, we may assume $\mathbf c$ very ample and $k$ algebraically closed. Then the base change to the function field $F$ of the generic curve does not change the heights (see \ref{reduction to generic curve}). By the usual specialization argument, the subvariety $Y \otimes_{\overline K} \overline F$ remains generic in $X \otimes_K F$ and hence \eqref{semipositive fundamental inequality for Y} follows  from the corresponding inequality over $F$. 

We apply \eqref{semipositive fundamental inequality for Y} now with $Y=Y_m$, $\overline L=(L, \canmetr)$  and we pass to the limit with respect to $m$. Setting $\ve=0$ and using non-negativity of the height, we deduce from \eqref{small for subvarieties} that $\hat{h}_L(X)=0$. For arbitrary $\ve>0$, we get
\begin{equation*}
\ve  \frac{h_{{\overline L}, \cdots ,{\overline L}, \overline O_X^f}(X)}{\deg_L(X)} \leq
\ve \liminf_m \frac{h_{{\overline L}, \dots ,{\overline L},\overline O_X^f}(Y_m)}{\deg_L(Y_m)} + O(\ve^2).
\end{equation*}
Letting $\ve \to 0$, we conclude
\begin{equation}  \label{last inequality}
\frac{h_{{\overline L}, \cdots ,{\overline L}, \overline O_X^f}(X)}{\deg_L(X)} \leq \liminf_m \frac{h_{{\overline L}, \dots ,{\overline L},\overline O_X^f}(Y_m)}{\deg_L(Y_m)}.
\end{equation}

Now let $f$ be a formal function on $\Xan$ inducing a formal $M_B$-metric $\metr_f$ on $O_X$ as at the beginning of this section. Then Theorem \ref{heights} shows
\begin{equation} \label{height of X in terms of chern form}
h_{\overline L, \dots, \overline L, \overline O_X^f}(X)=\mu(v)\int_{\Xan} f c_1(\overline L)^{\wedge d}. 
\end{equation}
Let us choose a finite normal subextension $K'/K$ of ${\overline K}/K$ such that $Y$ is defined over $K'$. Again Theorem \ref{heights} yields
\begin{equation} \label{height of Y in terms of chern form}
h_{\overline L, \dots, \overline L, \overline O_X^f}(Y_m)=\sum_{w|v} \mu(w)\int_{(Y_m)_w^{\rm an}} f c_1(\overline L)^{\wedge t}, 
\end{equation}
where $w$ ranges over all places of $K'$ over $v$. Since $K'/K$ is normal, ${\rm Aut}(K'/K)$ acts transitively on these places and hence $\mu(w)$ is independent of $w$. Moreover, $(Y_m)_w^{\rm an}$ ranges over the analytic spaces associated to the conjugates $Y_m^\sigma \in O(Y)$ and every $Y_m^\sigma$ is attained the same number of times. We conclude that
\begin{equation} \label{height and orbit}
h_{\overline L, \dots, \overline L, \overline O_X^f}(Y_m)=\lambda_m \sum_{Y_m^\sigma \in O(Y)} \int_{(Y_m^\sigma)_v^{\rm an}} f c_1(\overline L)^{\wedge t}
\end{equation}
for some $\lambda_m >0$. To determine the number $\lambda_m$, we use the constant function $f =1$ in \eqref{height of Y in terms of chern form}. Then the integral is equal to $\deg_L((Y_m)_w^{\rm an})=\deg_L(Y)$ by Proposition \ref{Chambert-Loir's measures for cycles}. The normalizations satisfy $\sum_{w|v}\mu(w)=\mu(v)$ (see \cite{BG}, Example 1.4.13) and hence we get $\lambda_m=\mu(v)|O(Y_m)|^{-1}$. Using \eqref{height of X in terms of chern form} and \eqref{height and orbit},  
Theorem \ref{equidistribution of small subvarieties} follows from \eqref{last inequality} in the same way as we proved Theorem \ref{main theorem}. \qed

\begin{rem} \rm \label{dynamical system for abelian variety}
If $W$ is an abelian variety $A$ over $K$ and $L$ is an ample even line bundle, then the theorem of the cube shows that multiplication by $m \in \ndop \setminus \{0,1\}$ leads to a dynamical system $(A,[m],L)$ with canonical height equal to the N\'eron--Tate height. For a $d$-dimensional irreducible closed subvariety $X$ of $A$,  the equidistribution measure $c_1(L|_X,\metr_{{\rm can},v})^{\wedge d}$ is explicitly described in terms of convex geometry by \cite{Gu6}, Theorem 6.7.
 \end{rem}
 
\begin{rem} \rm 
If we use the multiplicative torus $\Tor$ over $K$ instead of an abelian variety, then $\xb \mapsto \xb^m$ extends to a dynamical system $(\pdop_K^n,\Phi,O_{\pdop^n}(1))$ with canonical height equal to the classical Weil height of $\pdop_K^n$. Since $\pdop_K^n$ has good reduction with respect to any $v \in M_K$, the canonical equidistribution measure in Theorem \ref{equidistribution of small subvarieties} is the Dirac measure in the unique point of $(\pdop_K^n)_v^{\rm an}$ whose reduction modulo $v$ is the generic point. In the number field case, this non-archimedean analogue of Bilu's theorem (\cite{Bi}, Theorem 1.1) was proved by Chambert-Loir (see \cite{Ch}, Exemple 3.2).
\end{rem}

{\small Walter Gubler, Fachbereich Mathematik, Humboldt Universit\"at zu Berlin,
 Unter den Linden 6, D-10099 Berlin, walter.gubler@mathematik.uni-dortmund.de}

\begin{thebibliography}{EGA III}     
%\bibitem[Bei]{Bei}{A. Beilinson, {\it Height pairing between algebraic cycles}, Current trends in arithmetical algebraic geometry, Contemporary Mathematics Vol {\bf 67}, ed. K. Ribet, AMS 1985.}

\bibitem[Au]{Au}{P. Autissier: \'Equidistribution de sous-vari\'et\'es de petite hauteur. J. Th\'eor. Nombres Bordx. 18, No.1, 1--12 (2006).}

\bibitem[BI]{BI}{M. Baker, S.-I. Ih: Equidistribution of small subvarieties of an abelian variety. New York J. Math. 10, 279--289 (2004).}

\bibitem[Ber1]{Ber}{V.G. Berkovich: Spectral theory and analytic geometry over non-archimedean fields.  Mathematical Surveys and Monographs, 33. Providence, RI:  AMS (1990).}

%\bibitem[Ber2]{Ber2}{V.G. Berkovich:  \'Etale cohomology for non-archimedean analytic spaces. Publ. Math. IHES 78, 5--161 (1993).}

%\bibitem[Ber3]{Ber3}{V.G. Berkovich: Vanishing cycles for formal schemes. Invent. Math.  115, No.3, 539--571 (1994).}

%\bibitem[Ber4]{Ber4}{V.G. Berkovich:  Smooth $p$-adic analytic spaces are locally contractible. Invent. Math.  137, No.1, 1--84 (1999).}

%\bibitem[Ber5]{Ber5}{V.G. Berkovich: Smooth $p$-adic analytic spaces are locally contractible. II. Adolphson, Alan (ed.) et al., Geometric aspects of Dwork theory. Vol. I. Berlin: de Gruyter. 293--370 (2004).}

%\bibitem[BiGr]{BiGr}{R. Bieri, J.R.J. Groves: The geometry of the set of characters induced by valuations. J. Reine Angew. Math. 347, 168--195 (1984).}

\bibitem[Bi]{Bi}{Y. Bilu: Limit distribution of small points on algebraic tori. Duke Math. J. 89, No.3, 465--476 (1997).}



\bibitem[BG]{BG}{E. Bombieri, W.Gubler: Heights in Diophantine geometry. Cambridge University Press (2006).}

%\bibitem[Bl]{Bl}{S. Bloch, {\it Height pairings for algebraic cycles}, J. Pure Appl. Algebra {\bf 34} (1984), 119-145.}

%\bibitem[BlGS]{BlGS}{S. Bloch, H. Gillet, C. Soul\'e, {\it Non-archimedean Arakelov theory}, Publ. Math.IHES {\bf 78} (1995), 427-485.}

%\bibitem[BH]{BH}{T. Bloom, M. Herrera, {\it De Rham cohomology of an analytic space}, Invent. Math. {\bf 7} (1969), 275-296.}

%\bibitem[Bo]{Bo}{S. Bosch: Zur Kohomologietheorie rigid analytischer R\"aume. Manuscr. Math. 20, 1--27 (1977).}

\bibitem[BGR]{BGR}{S. Bosch, U. G\"untzer, R. Remmert: Non-Archimedean analysis. A systematic approach to rigid analytic geometry. Grundl. Math. Wiss., 261. Berlin etc.: Springer Verlag (1984).}

%\bibitem[BL1]{BL1}{S. Bosch, W. L\"utkebohmert: Stable reduction and uniformization of abelian varieties I. Math. Ann. 270, 349--379 (1985).}

%\bibitem[BL1]{BL1}{S. Bosch, W. L\"utkebohmert: N\'eron models from the rigid analytic viewpoint. J. Reine Angew. Math. 364, 69--84 (1986).}

%\bibitem[BL2]{BL2}{S. Bosch, W. L\"utkebohmert: Degenerating abelian varieties. Topology 30, No.4, 653--698 (1991).}

\bibitem[BL]{BL}{S. Bosch, W. L\"utkebohmert: Formal and rigid geometry. I: Rigid spaces. Math. Ann. 295, No.2, 291--317 (1993).}

%\bibitem[BL4]{BL4}{S. Bosch, W. L\"utkebohmert: Formal and rigid geometry. II: Flattening techniques. Math. Ann. 296, No.3, 403--429 (1993).}


%\bibitem[BoGS]{BoGS}{J.-B. Bost, H. Gillet, C. Soul\'e, {\it Heights of projective varieties and positive Green forms}, J. AMS {\bf 7/4} (1994), 903-1027.}

%\bibitem[Bou]{Bou}{N. Bourbaki: \'El\'ements de Math\'ematique. Fasc. XXX. Alg\`ebre commutative. Chap. 5: Entiers. Chap. 6: Valuations. Actualit\'es scientifiques et industrielles. 1308.  Paris: Hermann, 207 (1964).}

\bibitem[Ch]{Ch}{A. Chambert-Loir: Mesure et \'equidistribution sur les espaces de Ber\-ko\-vich. J. Reine Angew. Math. 595, 215--235 (2006).}

\bibitem[Fa]{Fa}{X.W.C. Faber: Equidistribution of dynamically small subvarieties over the function field of a curve. Preprint  at arXiv:math.NT:0801.4811v2.}

%\bibitem[CR]{CR}{T. Chinburg, R. Rumely, {\it The capacity pairing}, J. reine angew. Math. {\bf 434} (1993), 1-44.}

%\bibitem[EKL]{EKL}{M. Einsiedler, M. Kapranov, D. Lind: Non-archimedean amoebas and tropical varieties. J. Reine Angew. Math. 601, 139--157 (2006).}


%\bibitem[FC]{FC}{G. Faltings, C.-L. Chai: Degeneration of abelian varieties. Ergebnisse der  Mathematik und ihrer Grenzgebiete. 3. Folge, Bd. 22. Berlin etc.: Springer-Verlag (1989).}

%\bibitem[Fa1]{Fa1}{G. Faltings, {\it Calculus on arithm. surfaces}, Ann. of Math. {\bf 119} (1984), 387-424.}


%\bibitem[Fa2]{Fa2}{G. Faltings, {\it Diophantine approximation on abelian varieties}, Ann. Math. {\bf 133} (1991), 549-576.}

%\bibitem[FvdP]{FvdP}{J. Fresnel, M. van der Put: Rigid analytic geometry and its applications. Progress in Mathematics, 218. Boston, MA: Birkh\"auser. xii, 296 p. (2004).}

%\bibitem[Fu1]{Fu1}{W. Fulton: Intersection theory. Ergebnisse der  Mathematik und ihrer Grenzgebiete. 3. Folge, Bd. 2. Berlin etc.: Springer-Verlag (1984).}

%\bibitem[Fu2]{Fu2}{W. Fulton: Introduction to toric varieties. The 1989 W. H. Roever lectures in geometry. Annals of Math. Studies, 131.  Princeton, NJ: Pinceton University Press (1993).}

%\bibitem[GS1]{GS1}{H. Gillet, C. Soul\'e: Amplitude arithm\'etique. C.R. Acad Sci. Paris {\bf 307} (1988), 887-890.}


%\bibitem[GS2]{GS2}{H. Gillet, C. Soul\'e, {\it Arithmetic intersection theory}, Publ. Math. IHES {\bf 72} (1990), 94-174.}



%\bibitem[GH]{GH}{P. Griffiths, J. Harris, {\it Principles of algebraic geometry}, Pure and Applied Mathematics, Wiley, New York, 1978.}



%\bibitem[EGA II]{EGA II}{A. Grothendieck, J. Dieudonn\'e: \'El\'ements de g\'eometrie alg\'ebrique. II: \'Etude globale \'el\'ementaire de quelques classes de morphismes.  Publ. Math. IHES, 8, 1--222 (1961).}


%\bibitem[EGA IV]{EGA IV}{A. Grothendieck, J. Dieudonn\'e: \'El\'ements de g\'eometrie alg\'ebrique. IV: \'Etude locale des sch\'emas et des morphismes de sch\'emas (Quatrieme partie). Publ. Math. IHES, 1--361 (1967).}

%\bibitem[Gu1]{Gu1}{W. Gubler, {\it H\"ohentheorie} (mit einem Appendix von J\"urg Kramer), Math. Ann. {\bf 298} (1994), 427-455.}

%\bibitem[Gu1]{Gu1}{W. Gubler: Heights of subvarieties over $M$-fields. F. Catanese (ed.), Arithmetic geometry. Proceedings of a symposium, Cortona, 1994. Cambridge: Cambridge University Press.  Symp. Math.  37, 190-227 (1997).}

\bibitem[Gu1]{Gu2}{W. Gubler: Local heights of subvarieties over non-archimedean fields. J. Reine Angew. Math. 498, 61-113 (1998).} 

\bibitem[Gu2]{Gu3}{W. Gubler: Local and canonical heights of subvarieties. 
Ann. Sc. Norm. Super. Pisa, Cl. Sci. (Ser. V), 2, No.4., 711--760 (2003).}

\bibitem[Gu3]{Gu4}{W. Gubler: Tropical varieties for non-archimedean analytic spaces. Invent. Math. 169, No.2, 321--376 (2007).}

\bibitem[Gu4]{Gu5}{W. Gubler: The Bogomolov conjecture for totally degenerate abelian varieties. Invent. Math. 169, No.2, 377--400 (2007).}

\bibitem[Gu5]{Gu6}{W. Gubler: Non-archimedean canonical measures on abelian varieties. Preprint available at arXiv:math.NT:0801.4503v1.}

%\bibitem[Ha]{Ha}{U. Hartl: Semi-stable models for rigid-analytic spaces. Manuscr. Math. 110, No.3, 365--380 (2003).}

%\bibitem[HL]{HL}{U. Hartl, W. L\"utkebohmert: On rigid-analytic Picard varieties. J. Reine Angew. Math. 528, 101--148 (2000).} 

%\bibitem[Ha]{Ha}{R. Hartshorne, {\it Algebraic geometry}, GTM 52, Springer, Berlin, Heidelberg, New York, 1977.}

%\bibitem[Hr]{Hr}{P. Hriljac, {\it Heights and Arakelov's intersection theory}, Amer. J. Math. {\bf 167}/1. 23-38.}

%\bibitem[Ja]{Ja}{N. Jacobson, {\it Basic algebra II}, $2^{nd}$ edition, Freeman, New York 1989.}

%\bibitem[dJ]{dJ}{A. J. de Jong: Smoothness, semi-stability and alterations. Publ. Math. IHES 83, 51--93 (1996).}

%\bibitem[KKMS]{KKMS}{G. Kempf, F. Knudsen, D. Mumford, B. Saint-Donat: Toroidal embeddings. I. LNM 339. Berlin etc.: Springer-Verlag (1973).}

%\bibitem[Ki]{Ki}{J. King, {\it The currents defined by algebraic varieties}, Acta Math. {\bf 127} (1971), 185-220.}



%\bibitem[Kl]{Kl}{S. L. Kleiman: Toward a numerical theory of ampleness. Ann. Math. (2)  84, 293--344 (1966).}

%\bibitem[Ku1]{Ku1}{K. K\"unnemann: Projective regular models for abelian varieties, semistable reduction, and the height pairing. Duke Math. J. 95, No. 1, 161--212 (1998).}

%\bibitem[Ku2]{Ku2}{K. K\"unnemann: Height pairings for algebraic cycles on abelian varieties. Ann. Sci. \'Ec. Norm. Sup\'er. (4) 34, No. 4, 503--523 (2001).}


%\bibitem[Ko]{Ko}{U. K\"opf, {\it \"Uber eigentliche Familien algebraischer Variet\"aten \"uber affinoiden R\"aumen}, (German) Schr. Math. Inst. Univ. M\"unster.}

%\bibitem[La]{La}{S. Lang, {\it Fundamentals of diophantine geometry}, Berlin Heidelberg New York, Springer 1977.}

%\bibitem{La5} S. Lang, {\it Abelian Varieties.} Interscience Publishers, New York 1959. x+169 pp. Reprinted, Springer-Verlag, New York-Berlin 1983. xii+256 pp.

\bibitem[Lan]{Lan}{ S. Lang: Introduction to algebraic geometry. Interscience Publishers, New York (1958).}

\bibitem[Laz]{Laz}{ R. Lazarsfeld: Positivity in algebraic geometry. I. Classical setting: line bundles and linear series. Ergebnisse der Mathematik und ihrer Grenzgebiete. 3. Folge, Bd. 48. Berlin: Springer (2004).}

%\bibitem[McM]{McM}{P. McMullen: Duality, sections and projections of certain Euclidean tilings. Geom. Dedicata 49, No.2, 183--202 (1994).}

%\bibitem[Mi]{Mi}{G. Mikhalkin: Amoebas of algebraic varieties and tropical geometry. S. Donaldson (ed.) et al., Different faces of geometry. New York, NY: Kluwer Academic/Plenum Publishers. Int. Math. Ser. 3,  257--300 (2004).}



%\bibitem[Moi]{Moi}{B. G. Moishezon, {\it On $n$-dimensional compact varieties with $n$ algebraically independent meromorphic functions}, Amer. math. Soc. Translations {\bf 63} (1967), 51-177.}

%\bibitem[Mor]{Mor}{A. Moriwaki, {\it Arithmetic height functions over finitely generated fields}, Invent. Math. {\bf 140} (2000), 101-142.}

%\bibitem{Mu} D. Mumford, {\it Abelian Varieties.} Published for the Tata Institute of Fundamental Research Studies in Mathematics, No. 5. Oxford University Press, London 1970. viii+242 pp. 

%\bibitem[Mu]{Mu}{D. Mumford: An analytic construction of degenerating abelian varieties over complete rings. Compositio Math. 24, 239--272 (1972).}

%\bibitem[Ne]{Ne}{A. N\'eron, {\it Quasi-functions et hauteurs sur les vari\`etes ab\'eliennes}, Ann. of Math. {\bf 82} (1965), 249-331.}


%\bibitem[Ph1]{Ph1}{P. Philippon,{\it  Crit\`eres pour l'ind\'ependance alg\'ebrique}, Publ. Math. IHES {\bf 64} (1986), 5-52.}

%\bibitem[Ph2]{Ph2}{P. Philippon, {\it Sur des hauteurs alternatives I}, Math. Ann. {\bf 289} (1991), 255-283.}

%\bibitem[Ray]{Ray}{M. Raynaud, {\it Faisceaux amples sur les schemas en groupes et les espaces homog\'enes} Lecture Notes in Mathematics {\bf 119}, Springer-Verlag, 1970.}

%\bibitem[Oda]{Oda}{T. Oda: Convex bodies and algebraic geometry. An introduction to the theory of toric varieties. Ergebnisse der Mathematik und ihrer Grenzgebiete. 3. Folge, Bd. 15. Berlin etc.: Springer-Verlag. viii, 212 p. (1988).}

%\bibitem[Ru1]{Ru1}{W. Rudin: Functional Analysis. McGraw-Hill 
%Series in Higher Mathematics. New York etc.: McGraw-Hill. xiii, 397 p. (1973).}

%\bibitem[Ru2]{Ru2}{W. Rudin: Real and complex analysis.
%3rd. edition. New York, NY: McGraw-Hill. xiv, 416 p. (1987).}



\bibitem[SUZ]{SUZ} L. Szpiro, E. Ullmo and S. Zhang:  Equir\'epartition des petits points.  Invent.  Math.  127, 337--347 (1997).

\bibitem[Ullm]{Ullm}{E. Ullmo: Positivit\'e et discr\'etion des points alg\'ebriques des courbes. Ann. Math. (2) 147, No.1, 167--179 (1998).}

\bibitem[Ullr]{Ullr}{P. Ullrich:  The direct image theorem in formal and rigid geometry. Math. Ann. 301, No.1, 69--104 (1995).}

%\bibitem[Vo]{Vo}{P. Vojta, {\it Diophantine approximations and value distribution theory}, Springer lecture notes {\bf 1239} (1987).}

%\bibitem[We]{We}{A. Weil, {\it Arithmetic on algebraic varieties}, Ann. Math. {\bf 53} (1951), 412-444.}


%\bibitem[Zh1]{Zh1}{S. Zhang, {\it Admissible pairing on a curve}, Invent. Math. {\bf 112} (1993), 171-193.}

%\bibitem[Zh2]{Zh2}{S. Zhang, {\it Small points and adelic metrics}, J. Alg. Geometry {\bf 4} (1995), 281-300.}

%\bibitem[Zh3]{Zh3}{S. Zhang, {\it Equidistribution of small points on abelian varieties}, Ann. of Math. {\bf 147} (1998), 159-165.}

%\bibitem{Zh4} S. Zhang, Positive line bundles on arithmetic varieties, {\it J. Amer. Math. Soc.} {\bf 8} (1995), 187--221.

\bibitem[Yu]{Yu}{X. Yuan: Positive line bundles over arithmetic varieties. Preprint available at  arXiv:math.NT:0612424v1.}



\bibitem[Zh]{Zh}{S. Zhang: Equidistribution of small points on abelian varieties. Ann. Math. (2) 147, No.1, 159--165 (1998).}

\end{thebibliography}
\end{document}